\documentclass[10pt,a4paper]{article}

\usepackage{amsmath,amssymb,amsthm,mathtools}
\usepackage[T1]{fontenc}
\usepackage{lmodern}
\usepackage[
  textwidth=453.11pt,
  textheight=636.86pt,
  top=112.71pt,
  footskip=28.46pt
]{geometry}
\usepackage{empheq}
\usepackage{bm}
\usepackage{enumitem}
\usepackage{hyperref}
\usepackage{bookmark}
\usepackage{xcolor}
\usepackage{comment}

\let\articlenormalsize\normalsize
\renewcommand{\normalsize}{%
  \articlenormalsize
  \fontsize{10}{16}\selectfont
  \setlength{\abovedisplayskip}{11pt plus 3pt minus 4pt}%
  \setlength{\belowdisplayskip}{9pt plus 3pt minus 4pt}%
  \setlength{\abovedisplayshortskip}{0pt plus 3pt}%
  \setlength{\belowdisplayshortskip}{9pt plus 3pt minus 4pt}%
}
\let\articlesmall\small
\renewcommand{\small}{\articlesmall\fontsize{9}{13}\selectfont}
\let\articlefootnotesize\footnotesize
\renewcommand{\footnotesize}{\articlefootnotesize\fontsize{8}{11}\selectfont}
\normalsize
\everydisplay\expandafter{\the\everydisplay\fontsize{10}{12}\selectfont}

\hypersetup{
  colorlinks=true,
  linkcolor=blue,
  citecolor=blue,
  urlcolor=blue
}

\mathtoolsset{
  showonlyrefs=true,
  showmanualtags=true
}
\allowdisplaybreaks

\theoremstyle{plain}
\newtheorem{theorem}{Theorem}[section]
\newtheorem{proposition}[theorem]{Proposition}
\newtheorem{corollary}[theorem]{Corollary}
\newtheorem{lemma}[theorem]{Lemma}
\theoremstyle{definition}
\newtheorem{remark}[theorem]{Remark}

\newcounter{hypothesis}
\renewcommand{\thehypothesis}{H\arabic{hypothesis}}
\newenvironment{hypothesis}[1]{%
  \refstepcounter{hypothesis}%
  \par\medskip\noindent
  \textbf{(\thehypothesis) #1.}\quad
}{\par\medskip}
\newcommand{\hypref}[1]{(\ref{#1})}

\newcommand{\Td}{\mathbb T^d}

\newcommand{\dd}{\,dx}
\newcommand{\eps}{\varepsilon}
\newcommand{\norm}[2][]{\left\lVert #2\right\rVert_{#1}}
\newcommand{\abs}[1]{\left\lvert #1\right\rvert}

\newcommand{\baravg}[1]{\langle #1\rangle_{\bar m}}

\DeclareMathOperator{\Div}{div}

\title{Local turnpike property for a three-component mean field game system
arising from a resource harvesting model}
\author{%
  Mamoru Aihara\thanks{\raggedright Department of Mathematics, Graduate School of Science,
  Hokkaido University, Kita 10, Nishi 8, Kita-Ku, Sapporo, 060-0810, Hokkaido, Japan. 
  E-mail:
\href{mailto:aihara.mamoru.p7@elms.hokudai.ac.jp}{\texttt{aihara.mamoru.p7@elms.hokudai.ac.jp}}.}
  \and
  Shimpei Makida\thanks{\raggedright Faculty of Advanced Science and Technology,
  Ryukoku University, 1-5 Yokotani, Seta Oe-cho, Otsu, 520-2194, Shiga, Japan.
  E-mail: \href{mailto:makida@math.ryukoku.ac.jp}{\texttt{makida@math.ryukoku.ac.jp}}.}%
}
\date{}

\begin{document}
\maketitle

\begin{abstract}
We study a three-component mean field game system arising from a resource harvesting model.
It couples Hamilton--Jacobi and Fokker--Planck equations for the harvesters
with a reaction-diffusion equation for the resource.
The coupling is nonlocal in space and time.
For standard second-order mean field game systems, exponential turnpike estimates are known
under assumptions such as monotonicity of the coupling.
For the present system, Kobeissi, Mazari-Fouquer, and Ruiz-Balet proved algebraic convergence
to a stationary solution in a space--time averaged sense,
but no exponential turnpike estimate has been obtained.
We prove such an estimate near a stationary solution with strictly positive density components.
Inspired by the relative entropy method, we add the relative entropy of the resource density
to the standard energy for local turnpike estimates.
In the time derivative of this modified energy, the nonlocal coupling terms cancel.
The two key assumptions are coercivity of a quadratic form determined by the stationary resource density
and smallness of the relative density perturbations in $L^\infty$.
Under these assumptions, we bound the absolute value of the modified energy by two terms
that decay exponentially away from the initial and terminal times.
We derive weighted $L^2$ exponential turnpike estimates for the harvester density,
the resource density, and the gradient of the value function.
We give sufficient conditions for coercivity.
For each fixed $T>0$, classical solutions satisfying these assumptions exist
for initial and terminal data sufficiently close to the stationary solution.
We extend the modified energy estimate to a more general reaction term under the corresponding assumptions.
\end{abstract}
\section{Introduction}

Mean field game (MFG) theory studies noncooperative games with a large population of identical players.
The players interact through the distribution of the population.
The theory describes Nash equilibria in the limit as the number of players tends to infinity.
It was introduced independently by Lasry--Lions \cite{LL07}
and Huang--Malham\'{e}--Caines \cite{HMC06} and has since developed rapidly.
In basic models where each player's state is perturbed by Brownian motion,
the resulting MFG system consists of a Hamilton--Jacobi equation for the value function $u$ of a representative player
and a Fokker--Planck equation for the density $m$ of the players.
The first equation runs backward in time from a terminal condition,
while the second runs forward from an initial condition.
This forward-backward structure makes it difficult to study MFG systems
by the methods used for initial-value problems for standard parabolic systems.

On a long time interval $[0,T]$, solutions may deviate from a stationary solution
near the initial and terminal times because they must satisfy the initial and terminal conditions.
These short intervals near the endpoints are called temporal boundary layers.
Following terminology from economics and optimal control,
the turnpike property means that solutions remain close to a stationary solution
for most of the time interval, outside these boundary layers.
We say that the exponential turnpike property holds if the deviation is bounded by
$C(e^{-\sigma t}+e^{-\sigma(T-t)})$ for constants $C,\sigma>0$ independent of $T$.

Cardaliaguet--Lasry--Lions--Porretta \cite{CLLP12,CLLP13} studied the long-time behavior of MFG systems.
They studied convergence of long-time averages to stationary solutions of the ergodic system
under Lasry--Lions monotonicity of the coupling.
They also studied exponential estimates under additional assumptions
on the strength of monotonicity and the regularity of the coupling.
Exponential turnpike properties have also been studied under weaker forms of Lasry--Lions monotonicity
\cite{CP21,CCDE24}.
In contrast, for systems with a quadratic Hamiltonian, Cirant--De Bernardi \cite{CN26}
proved an exponential turnpike property for solutions near a stationary solution,
under a second-order strict positivity condition at that stationary solution
in place of global monotonicity.
This property is called the local turnpike property.
Their condition requires a quadratic form to be coercive on density perturbations with zero total mass.
This quadratic form contains the derivative of the coupling with respect to the density at the stationary solution.
Thus, the condition incorporates information about the coupling.
For a nonmonotone MFG of Kuramoto type with sufficiently strong interaction,
\cite{CC24} obtained local exponential estimates at each time
for the densities of reflection-symmetric solutions that remain near a self-organized stationary solution.
Very recently, Ghattassi--Kong--Masmoudi \cite{GKM26} extended local turnpike analysis
to two-population second-order MFG systems with quadratic Hamiltonians under a matrix-valued
coercivity condition accounting for cross-population interactions.
Their system consists of two coupled Hamilton--Jacobi--Fokker--Planck pairs, whereas the present
work concerns one such pair coupled to a resource reaction--diffusion equation.

Recent MFG models in economics and mathematical ecology include multicomponent systems
with additional unknowns or several populations.
Examples include price formation models in which prices are determined endogenously
through a market-clearing condition \cite{GS21,FT22,Fuj23}.
Another example is a model for two populations of large and small beetles,
with the probability of expression of a gene that determines body size as an additional unknown \cite{JLWZZ25}.
Kobeissi--Mazari-Fouquer--Ruiz-Balet \cite{KMR24} studied a harvesting model
that couples an MFG system with a reaction-diffusion equation for the resource population.
They interpreted the tragedy of the commons (resource depletion caused by competition) as a reversal of traveling waves in the reaction-diffusion equation.

In this paper, we study the MFG system introduced by
Kobeissi--Mazari-Fouquer--Ruiz-Balet \cite{KMR25} for a resource harvesting problem.
The spatial domain is the $d$-dimensional torus $\Td=\mathbb R^d/\mathbb Z^d$,
with periodic boundary conditions.
Let $\nu,\mu,\varepsilon>0$ be constants.
Let $K:\Td\to\mathbb R$ be the environmental function representing the spatially varying intrinsic growth rate
of the resource, and let $\rho:\Td\to\mathbb R$ be a convolution kernel.
We define convolution by
$
  (\rho*\theta)(x):=\int_{\Td}\rho(x-y)\theta(y)\,dy
$.
We consider the system (KMR)
\begin{align}
-\partial_t u-\nu\Delta u-\frac12|\nabla u|^2
&=\rho*\theta
&&\quad\text{in }(0,T)\times\Td,
\tag{\textup{KMR-1}}\label{eq:KMR-HJ}
\\
\partial_t m-\nu\Delta m+\nabla\cdot(m\nabla u)
&=0
&&\quad\text{in }(0,T)\times\Td,
\tag{\textup{KMR-2}}\label{eq:KMR-FP}
\\
\partial_t\theta-\mu\Delta\theta
&=\theta(K-\theta)-\varepsilon(\rho*m)\theta
&&\quad\text{in }(0,T)\times\Td
\tag{\textup{KMR-3}}\label{eq:KMR-theta}
\end{align}
with the initial and terminal conditions
\[
  u(T,\cdot)=u_T,\qquad
  m(0,\cdot)=m_0,\qquad
  \theta(0,\cdot)=\theta_0
  \quad\text{on }\Td.
\]
The players who harvest the resource are called harvesters.
Here, $u$ is the value function of a representative harvester, representing the maximal expected payoff,
$m$ is the harvester density, and $\theta$ is the resource density.
The constants $\nu$ and $\mu$ are the diffusion coefficients of the harvesters and the resource, respectively,
and $\varepsilon$ is the harvesting intensity.
The term $\rho*\theta$ on the right-hand side of \eqref{eq:KMR-HJ}
represents the resource available to harvesters at each location and gives their running payoff.
The term $\varepsilon(\rho*m)\theta$ in \eqref{eq:KMR-theta} describes resource depletion due to harvesting.

In (KMR), the distribution of harvesters affects the evolution of the resource.
The resource density, in turn, affects the harvesters' optimal behavior through the payoff $\rho*\theta$.
In particular, solving \eqref{eq:KMR-theta} expresses the resource density $\theta$
as a functional of the harvester distribution $m$.
The value $\theta(t)$ depends on the entire history of $m$ over $[0,t]$.
Thus, the coupling term $\rho*\theta$ in \eqref{eq:KMR-HJ} is nonlocal in both space and time,
whereas standard nonlocal couplings depend only on the distribution at the current time \cite{KMR25}.

The ergodic (stationary) system (sKMR) associated with (KMR) is
\begin{align}
\lambda-\nu\Delta\bar u-\frac12|\nabla\bar u|^2
&=\rho*\bar\theta
&&\quad\text{in }\Td,
\tag{\textup{sKMR-1}}\label{eq:sKMR-HJ}
\\
-\nu\Delta\bar m+\nabla\cdot(\bar m\nabla\bar u)
&=0
&&\quad\text{in }\Td,
\tag{\textup{sKMR-2}}\label{eq:sKMR-FP}
\\
-\mu\Delta\bar\theta
&=\bar\theta(K-\bar\theta)-\varepsilon(\rho*\bar m)\bar\theta
&&\quad\text{in }\Td
\tag{\textup{sKMR-3}}.\label{eq:sKMR-theta}
\end{align}
We call a tuple $(\bar u,\bar m,\bar\theta,\lambda)$ a stationary solution if it satisfies this system and
\[
  \int_{\Td}\bar u\dd=0,\qquad
  \bar m\geq0,\qquad\int_{\Td}\bar m\dd=1,\qquad
  \bar\theta\geq0,\qquad\bar\theta\not\equiv0.
\]
Here, the unknown constant $\lambda\in\mathbb R$ is the ergodic constant, which represents the long-time average payoff per unit time.
The condition $\int_{\Td}\bar u\dd=0$ fixes the additive constant in $\bar u$.
For each stationary solution, $(\bar u+\lambda(T-t),\bar m,\bar\theta)$ satisfies the equations of the time-dependent system (KMR).
We call this tuple the reference solution and compare solutions of (KMR) with it.

We impose the following conditions on the environmental function and the convolution kernel.
\begin{hypothesis}{Regularity of the environmental function and positivity of its integration}\label{hyp:env}
  For some $\beta\in(0,1)$, we have
  \[
    K\in C^{2+\beta}(\Td),\qquad\int_{\Td}K\dd>0.
  \]
\end{hypothesis}
\begin{hypothesis}{Regularity, nonnegativity, and symmetry of the convolution kernel}\label{hyp:even-kernel}
  The kernel $\rho\in C^\infty(\Td)$ is nonnegative and satisfies
  \[
    \rho(x)=\rho(-x),\qquad\int_{\Td}\rho\dd=1.
  \]
  In particular, $\rho$ is even.
\end{hypothesis}
Condition \hypref{hyp:env} allows $K$ to be nonpositive at some points.
Thus, regions $\{K\leq0\}$ unfavorable to resource growth may be present.

Existence of solutions to (sKMR) is studied in \cite[Theorem~1]{KMR25}.
We consider stationary solutions with positive resource density and use the following sufficient condition for their existence.
\begin{theorem}[Existence of a stationary solution with positive resource density]\label{thm_existence_ergo}
  Assume \hypref{hyp:env} and \hypref{hyp:even-kernel}, and let
  $0<\varepsilon<\int_{\Td}K\dd$.
  Then (sKMR) admits at least one classical solution $(\bar u,\bar m,\bar\theta,\lambda)$
  with $\bar m>0$ and $\bar\theta>0$.
\end{theorem}

The condition on $\varepsilon$ is sufficient for the existence of a stationary solution with positive resource density.
Indeed, for any probability density $\widehat m$, consider $K-\varepsilon(\rho*\widehat m)$.
By \hypref{hyp:even-kernel}, its spatial average is
\[
\int_{\mathbb T^d}Kdx-\varepsilon>0.
\]
Hence, the Rayleigh quotient estimate shows that the principal eigenvalue $\lambda_1(\widehat m)$
 of the operator $-\mu\Delta-K+\varepsilon(\rho*\widehat m)$ satisfies
\[
\lambda_1(\widehat m)
\le
\varepsilon-\int_{\mathbb T^d}K\ dx<0.
\]
The negativity of the principal eigenvalue allows us to apply \cite[Theorems 2.1(1) and 2.4]{BHR05}.
For the fixed density $\widehat m$, equation \eqref{eq:sKMR-theta} therefore has a unique positive solution $\theta_{\widehat m}$.

Next, using this resource density $\theta_{\widehat m}$, we solve \eqref{eq:sKMR-HJ} and \eqref{eq:sKMR-FP}
(see \cite[Section~2, p.~3560]{CLLP13} for the stationary two-component system).
Let $\mathcal T$ be the map that returns the resulting density component.
Schauder's fixed point theorem yields a fixed point $\bar m=\mathcal T(\bar m)$, and hence a stationary solution.

Our aim is to derive exponential turnpike estimates for solutions of (KMR) that remain,
throughout the time interval, near a stationary solution whose densities $\bar m,\bar\theta$ are both positive.
For the long-time behavior of this system, \cite[Theorem~4]{KMR25} establishes an $O(T^{-1})$ bound
on the squared $L^2((0,1)\times\Td)$ distance between the rescaled resource density $\theta(Ts,\cdot)$
and the stationary resource density, where $t=Ts$ and $s\in[0,1]$.
This result assumes that the initial resource density is the stationary resource density in the absence of harvesting,
that the harvesting intensity $\varepsilon$ is sufficiently small,
and that either the spatial dimension is one or the spatial variation of the environmental function $K$ is sufficiently small.
It gives a global stability result and also applies when the harvesters have no diffusion, $\nu=0$.
The estimate is averaged in time and does not give a decay rate at each time.
In contrast, for $\nu>0$, we assume coercivity of a quadratic form determined by the stationary resource density
(denoted below by $Q_\theta$), together with smallness of the relative perturbations
$(m-\bar m)/\bar m$ and $(\theta-\bar\theta)/\bar\theta$ in $L^\infty$.
We prove that the squared weighted $L^2$ norms of $m-\bar m$, $\theta-\bar\theta$, and $\nabla u-\nabla\bar u$
are bounded at each $t\in[0,T]$ by a linear combination of $e^{-\sigma t}$ and $e^{-\sigma(T-t)}$.
This has the same form as the exponential turnpike estimates for standard second-order MFG systems with monotone couplings
\cite{CLLP12,CLLP13}, and gives the decay rate naturally expected for the present system.
The results of \cite{KMR25} and the present paper are complementary, with the former covering a broader range of solutions
and the latter providing stronger estimates.

For the analysis below, we make the following assumptions on the classical solutions and their mass.
\begin{hypothesis}{Classical solutions and positivity of the stationary densities}\label{hyp:regularity}
  The tuples $(u,m,\theta)$ and $(\bar u,\bar m,\bar\theta,\lambda)$ are classical solutions
  of (KMR) and (sKMR), respectively.
  There are constants $C_m,C_\theta\geq1$ such that
  \[
    C_m^{-1}\leq\bar m(x)\leq C_m,\qquad
    C_\theta^{-1}\leq\bar\theta(x)\leq C_\theta
    \quad\text{on }\Td.
  \]
\end{hypothesis}
\begin{hypothesis}{Mass normalization}\label{hyp:mass}
  The initial density $m_0$ is nonnegative and satisfies $\int_{\Td}m_0\dd=1$.
\end{hypothesis}
We define the perturbation variables relative to the reference solution by
\[
  v:=u-\bar u-\lambda(T-t),\qquad
  n:=m-\bar m,\qquad z:=\theta-\bar\theta,
\]
\[
  r:=\frac{n}{\bar m},\qquad
  q:=\frac{z}{\bar\theta},\qquad
  A:=\nabla v-\nu\nabla r.
\]
Here, $r$ and $q$ are the relative perturbations of the stationary densities $\bar m$ and $\bar\theta$, respectively.
By \hypref{hyp:mass} and mass conservation, $\int_{\Td}n\dd=\int_{\Td}\bar m r\dd=0$.

\medskip\noindent\textbf{Main ideas and results.}
We now describe the construction of the modified energy $E$ and state the main results obtained from it.

The energy method starts from the cross energy arising from the duality between the Hamilton--Jacobi and Fokker--Planck equations \cite{CLLP12},
\[
  E_1(t):=-\int_{\Td}nv\dd.
\]
This corresponds to the quantity $\int_{\Td}(u_1-u_2)(m_1-m_2)\dd$ in the Lasry--Lions uniqueness argument.
Using the equations for the perturbations $(v,n,z)$ (see Section~2),
the zero-flux identity $\nu\nabla\bar m=\bar m\nabla\bar u$ (Lemma~\ref{lem:fluz}), and the evenness of $\rho$, we obtain
\[
  E_1'(t)
  =-\int_{\Td}\bar m|\nabla v|^2\dd
   -\frac12\int_{\Td}n|\nabla v|^2\dd
   +\int_{\Td}(\rho*n)z\dd.
\]
For standard MFG systems, Lasry--Lions monotonicity determines the sign of the coupling contribution,
which can then be treated as a dissipation term.
In the present system, the coupling is nonlocal in time, so monotonicity at each time generally fails.
Thus, the term $\int_{\Td}(\rho*n)z\dd$ cannot be treated in the same way.
The monotonicity property integrated over space and time in \cite[Proposition~1]{KMR25}
compares two resource densities with the same initial resource density.
It therefore does not apply to the comparison between a time-dependent solution and a stationary solution.

Following the modified cross energy construction in \cite{CN26}, we first add
$\nu\|n/\sqrt{\bar m}\|_{L^2(\Td)}^2$ to the cross energy and define
\[
  E_2(t):=-\int_{\Td}nv\dd
           +\nu\int_{\Td}\frac{n^2}{\bar m}\dd
         =-\int_{\Td}\bar m r(v-\nu r)\dd.
\]
This correction yields the dissipation term $-\nu^2\int_{\Td}\bar m|\nabla r|^2\dd$
involving the gradient of $r$, which was absent from $E_1'$. We obtain
\[
\begin{aligned}
  E_2'(t)
  ={}&-\int_{\Td}\bar m|\nabla v-\nu\nabla r|^2\dd
      -\nu^2\int_{\Td}\bar m|\nabla r|^2\dd
      +\int_{\Td}(\rho*n)z\dd+R_m(t),
  \\
  R_m(t):={}&-\frac12\int_{\Td}n|\nabla v|^2\dd
              +2\nu\int_{\Td}n\nabla r\cdot\nabla v\dd.
\end{aligned}
\]
The order of a term in the perturbations is the number of factors in the term
that involve $v,r,q$, their derivatives, or their convolutions.
Since $n=\bar m r$ and $z=\bar\theta q$, both dissipation terms above and the coupling term $\int_{\Td}(\rho*n)z\dd$ are quadratic.
The coupling term therefore cannot be treated as a higher-order error.
Moreover, the dissipation terms in $E_2'$ do not involve the resource perturbation $z$ and cannot absorb this coupling term.
We use the structure of the resource equation to cancel it.

To this end, we introduce the modified energy
\[
  E(t):=E_2(t)+\frac1\varepsilon H_\theta(t).
\]
That is, $E$ is obtained by adding to $E_2$ the relative entropy of the resource density
with weight $1/\varepsilon$, which is determined by the harvesting intensity,
\[
  H_\theta(t):=\int_{\Td}\bar\theta\{q-\log(1+q)\}\dd.
\]
Here, $H_\theta$ is the relative entropy of $\theta$ with respect to $\bar\theta$, defined for $\theta>0$.
It has the same form as the Lyapunov functionals used in the stability analysis of logistic equations and Lotka--Volterra systems.
The resource perturbation equation gives
\[
\begin{aligned}
  \frac1\varepsilon H_\theta'(t)
  ={}&-\frac\mu\varepsilon
       \int_{\Td}\bar\theta\frac{|\nabla q|^2}{(1+q)^2}\dd
       +\frac\mu\varepsilon
       \int_{\Td}\frac{q}{1+q}\nabla\bar\theta\cdot\nabla q\dd
  \\
     &-\frac1\varepsilon\int_{\Td}\bar\theta^2q^2\dd
       -\int_{\Td}z(\rho*n)\dd.
\end{aligned}
\]
The last term on the right-hand side exactly cancels the coupling term $\int_{\Td}(\rho*n)z\dd$ in $E_2'$.
After this cancellation, we obtain
\[
\begin{aligned}
  E'(t)={}&-\int_{\Td}\bar m|\nabla v-\nu\nabla r|^2\dd
            -\nu^2\int_{\Td}\bar m|\nabla r|^2\dd
  \\
          &-\frac\mu\varepsilon
            \int_{\Td}\bar\theta\frac{|\nabla q|^2}{(1+q)^2}\dd
            +\frac\mu\varepsilon
            \int_{\Td}\frac q{1+q}\nabla\bar\theta\cdot\nabla q\dd
  \\
          &-\frac1\varepsilon\int_{\Td}\bar\theta^2q^2\dd
            -\frac12\int_{\Td}n|\nabla v|^2\dd
            +2\nu\int_{\Td}n\nabla r\cdot\nabla v\dd.
\end{aligned}
\]
The last two terms on the right-hand side sum to $R_m(t)$.

Near the stationary solution, terms of order three or higher are small compared with the quadratic terms.
We therefore decompose $E'$ into a leading quadratic part and a remainder of order at least three.
The first two terms and the $q^2$ term on the right-hand side of $E'$ are already quadratic, while $R_m$ is cubic.
In the two terms involving the resource gradient, expand $(1+q)^{-2}$ and $q/(1+q)$ around $q=0$.
Their lowest-order terms, $1$ and $q$, yield quadratic terms, while the remaining terms have order at least three.
To organize these terms, define the flux dissipation
\[
  D_{\mathrm{flux}}(t):=\int_{\Td}\bar m|A|^2\dd
\]
and the quantity determined by the stationary resource density,
\[
  Q_\theta(\varphi):=
  \mu\int_{\Td}\bar\theta|\nabla\varphi|^2\dd
  -\mu\int_{\Td}\varphi\nabla\bar\theta\cdot\nabla\varphi\dd
  +\int_{\Td}\bar\theta^2\varphi^2\dd,
  \qquad\varphi\in H^1(\Td).
\]
The vector field $\bar mA=\bar m(\nabla v-\nu\nabla r)$ is the density flux in the Fokker--Planck equation
linearized around the stationary solution. Thus, $D_{\mathrm{flux}}$ is the squared weighted $L^2$ norm of this flux, with weight $1/\bar m$.
The quantity $-Q_\theta(q)/\varepsilon$ is the sum of the quadratic parts of the three resource terms.
The coefficients of $Q_\theta$ involve only $\mu$ and the stationary resource density $\bar\theta$.
In contrast to the quadratic form in the positivity condition of \cite{CN26}, it contains no density derivative of the coupling.
The second term has no fixed sign.
Thus, $Q_\theta$ is a quadratic form associated with the stationary solution, and it is not necessarily nonnegative.
Combining the density-gradient term with the resource contribution, we define the quadratic form
\[
  \mathcal{Q}(r,q):=
  \nu^2\int_{\Td}\bar m|\nabla r|^2\dd
  +\frac1\varepsilon Q_\theta(q).
\]
The quadratic part of $E'$ is exactly $-D_{\mathrm{flux}}-\mathcal{Q}$, so the energy identity becomes
\[
  E'(t)=-D_{\mathrm{flux}}(t)-\mathcal{Q}(r,q)+R_m(t)+R_\theta(t).
\]
Here, $R_\theta$ is the remainder after removing the quadratic parts of the two resource gradient terms,
\[
\begin{aligned}
  R_\theta(t):={}\frac\mu\varepsilon
    \int_{\Td}\bar\theta\left\{1-\frac1{(1+q)^2}\right\}
    |\nabla q|^2\dd
  +\frac\mu\varepsilon
    \int_{\Td}\left\{\frac q{1+q}-q\right\}
    \nabla\bar\theta\cdot\nabla q\dd.
\end{aligned}
\]
The leading quadratic part $D_{\mathrm{flux}}+\mathcal{Q}$ therefore represents the dissipation of $E$
for the system obtained by linearizing (KMR) around the stationary solution.
The terms $R_m$ and $R_\theta$ are higher-order remainders.
To control these remainders by the leading terms, we impose the following two conditions.
\begin{hypothesis}{Coercivity of $Q_\theta$}\label{hyp:spectral-gap}
  There exists $\lambda_\theta>0$ such that, for every $\varphi\in H^1(\Td)$,
  \[
    Q_\theta(\varphi)\geq\lambda_\theta
      \int_{\Td}\bar\theta\varphi^2\dd.
  \]
\end{hypothesis}
\begin{hypothesis}
{Smallness of the relative perturbations in \(L^\infty\)}\label{hyp:smallness}
Define the constant
\[
C_V
:=
\frac32\|\bar\theta\|_{L^\infty}
+
\frac12\|K\|_{L^\infty}
+
\frac{\varepsilon}{2}\|\rho\|_{L^\infty}.
\]
It depends on the $L^\infty$ norms of $K$, $\bar \theta$, and $\rho$, and on the harvesting intensity $\varepsilon$.
Set
\[
C_R:=
\max\left\{2, \left(\frac{\mu}{\lambda_\theta}\left\|\frac{\nabla \bar \theta}{\bar \theta}\right\|_{L^\infty} ^2+7(1+\frac{C_V}{\lambda_{\theta}})\right) \right\}.
\]
There exists a constant \(\delta\in(0,1/2]\) such that
\[
C_R \,\delta \le \frac{1}{2},
\qquad
\sup_{t\in[0,T]}
\left(
\|r(t)\|_{L^\infty(\mathbb T^d)}
+\|q(t)\|_{L^\infty(\mathbb T^d)}
\right)
\le \delta.
\]
\end{hypothesis}
Condition \hypref{hyp:spectral-gap} requires the quadratic form $Q_\theta$ to be coercive
with respect to the weighted $L^2$ norm $\|\sqrt{\bar\theta}\,\varphi\|_{L^2(\Td)}$.
We refer to it as the coercivity condition and discuss it further at the end of the introduction.

Under \hypref{hyp:spectral-gap} and \hypref{hyp:smallness}, the remainders satisfy
\[
  |R_m(t)|+|R_\theta(t)|
  \leq C_R\delta\{D_{\mathrm{flux}}(t)+\mathcal{Q}(r(t),q(t))\}
\]
with the constant $C_R$ defined in \hypref{hyp:smallness}, which is independent of $T$ (see Section~2).
The smallness of $\delta$ required in \hypref{hyp:smallness} allows us to absorb the remainders into the leading terms.
This gives a dissipation estimate at each time without a monotonicity assumption on the coupling.
\begin{proposition}[Dissipation estimate for the modified energy]\label{prop:local-energy-di}
  Assume \hypref{hyp:even-kernel} and \hypref{hyp:regularity}--\hypref{hyp:smallness}.
  Then there exists a constant $C_E>0$, independent of $T$, such that
  \begin{equation}
    -E'(t)\geq \frac{1}{2}\{D_{\mathrm{flux}}(t)+\mathcal{Q}(r(t),q(t))\},
    \label{eq:main-dissipation}
  \end{equation}
  \begin{equation}
    |E(t)|\leq C_E\{D_{\mathrm{flux}}(t)+\mathcal{Q}(r(t),q(t))\}.
    \label{eq:main-energy-control}
  \end{equation}
  Consequently, setting $\sigma:=\frac{1}{2C_E}>0$, we obtain
  \begin{equation}
    E'(t)\leq-\sigma|E(t)|.
    \label{eq:main-signed-energy}
  \end{equation}
\end{proposition}
The estimates of Proposition~\ref{prop:local-energy-di} also hold when the reaction term $\theta(K-\theta)$
in the resource equation is replaced by a general function $G(x,\theta)$,
provided that the reaction term has suitable regularity and that \hypref{hyp:spectral-gap} and \hypref{hyp:smallness}
are modified accordingly (Proposition~\ref{prop:general-reaction}).

Integrating \eqref{eq:main-signed-energy} over $[0,t]$ when $E(t)\geq0$
and over $[t,T]$ when $E(t)<0$, we obtain
\[
  |E(t)|\leq |E(0)|e^{-\sigma t}+|E(T)|e^{-\sigma(T-t)}.
\]
Combining this bound with differential inequalities for the weighted $L^2$ norms of the components yields the following result.
\begin{proposition}[Weighted $L^2$ turnpike estimates]\label{prop:L2-turnpike_main}
  Assume \hypref{hyp:even-kernel} and \hypref{hyp:regularity}--\hypref{hyp:smallness}.
  Then there exist constants $C,\sigma>0$, independent of $T$, such that, for every $t\in[0,T]$,
  \[
    \left\|\frac{m(t)-\bar m}{\sqrt{\bar m}}\right\|_{L^2(\Td)}^2
    \leq\lambda_1e^{-\sigma t}+\lambda_2e^{-\sigma(T-t)},
  \]
  \[
    \left\|\frac{\theta(t)-\bar\theta}{\sqrt{\bar\theta}}\right\|_{L^2(\Td)}^2
    \leq\lambda_3e^{-\sigma t}+\lambda_4e^{-\sigma(T-t)},
  \]
  \[
    \|\sqrt{\bar m}\,\nabla v(t)\|_{L^2(\Td)}^2
    \leq\lambda_5e^{-\sigma t}+\lambda_6e^{-\sigma(T-t)}.
  \]
  The coefficients may be chosen as
  \[
  \begin{alignedat}{2}
    \lambda_1&=C\left(|E(0)|+\|\sqrt{\bar m}\,r(0)\|_{L^2}^2\right),\qquad
    &\lambda_2&=C|E(T)|,\\
    \lambda_3&=C\left(|E(0)|+\|z(0)/\sqrt{\bar\theta}\|_{L^2}^2\right),\qquad
    &\lambda_4&=C|E(T)|,\\
    \lambda_5&=C|E(0)|,\qquad
    &\lambda_6&=C\left(|E(T)|+\|\sqrt{\bar m}\,\nabla v(T)\|_{L^2}^2\right).
  \end{alignedat}
  \]
\end{proposition}
The strategy of deriving the turnpike property from an energy estimate is shared with
Cardaliaguet--Lasry--Lions--Porretta \cite{CLLP12}.
In particular, our derivation of weighted $L^2$ estimates for each component from the modified energy dissipation estimate
follows Cirant--De Bernardi \cite{CN26}.
The $L^2$ turnpike estimates in \cite{CN26} assume uniform smallness of $\nabla v$ in $T$.
Here, the cubic remainder terms involving $\nabla v$ in the energy derivative are absorbed
using the smallness of the relative density perturbation $r$.
For the turnpike estimate on $\nabla v$, we use the time-dependent density $m$ as a weight.
This cancels the nonlinear transport terms arising from the Fokker--Planck and Hamilton--Jacobi equations.
Thus, Proposition~\ref{prop:L2-turnpike_main} gives $L^2$ turnpike estimates
without assuming an a priori $L^\infty$ bound on $\nabla v$ that is uniform in $T$.

The endpoint values $E(0)$ and $E(T)$ in Proposition~\ref{prop:L2-turnpike_main}
also depend on unknown endpoint values of the solution.
If we additionally assume $\|\nabla v(0)\|_{L^\infty(\mathbb{T}^d)}\le 1$,
their absolute values can be bounded in terms of the prescribed initial and terminal data $m_0,\theta_0,u_T$
and the stationary solution.
The coefficients in the turnpike estimates can then be controlled by the deviations of the initial and terminal data
from the stationary solution (see Lemmas~\ref{lem:E-estimate} and \ref{lem:L2-turnpike-}).

Next, we study sufficient conditions for \hypref{hyp:spectral-gap}.
Using the stationary resource equation \eqref{eq:sKMR-theta} and the classical maximum principle,
we obtain a condition for coercivity of $Q_\theta$ in terms of the lower and upper bounds
of the effective growth rate $F:=K-\varepsilon(\rho*\bar m)$.
This is the environmental function minus the depletion rate due to harvesting at the stationary solution.
In particular, if $K\equiv K_0>0$ and $0<\varepsilon<K_0$, the identity $\rho*1=1$ shows directly
that the spatially homogeneous tuple $(\bar u,\bar m,\bar\theta,\lambda)=(0,1,K_0-\varepsilon,K_0-\varepsilon)$
is a stationary solution.
Since $\nabla\bar\theta=0$, condition \hypref{hyp:spectral-gap} holds with $\lambda_\theta=K_0-\varepsilon$.
This sufficient condition is consistent with the resource dynamics in a homogeneous environment without harvesting
($\varepsilon=0$, $K\equiv K_0$).
Indeed, when $\varepsilon=0$ and $K\equiv K_0$, solutions of the initial-value problem for the Fisher equation
with positive initial data are known to converge exponentially to the steady state.

We also prove the existence of solutions near the stationary solution that satisfy
\hypref{hyp:regularity} and \hypref{hyp:smallness}, as assumed in Propositions~\ref{prop:local-energy-di} and \ref{prop:L2-turnpike_main}.
The proof uses the Leray--Schauder fixed point theorem.
We construct a homotopy map for the system of perturbations around the stationary solution.
An improvement of the $L^\infty$ turnpike estimate (Proposition~\ref{prop:n-Linfty-improvement})
provides the a priori estimates needed to obtain a fixed point.
This argument is based on the fixed point method of \cite{CN26}.
For the present system, positivity, regularity, and compactness estimates are also needed for the resource reaction-diffusion equation.
\begin{proposition}[Existence of solutions near a stationary solution]\label{prop_exsist}
  Assume \hypref{hyp:env} and \hypref{hyp:even-kernel}.
  Let $(\bar u,\bar m,\bar\theta,\lambda)$ be a classical solution of (sKMR)
  satisfying $\bar m>0$, $\bar\theta>0$, and \hypref{hyp:spectral-gap}.
  Suppose that the initial and terminal data belong to $C^{2+\beta}(\Td)$,
  with $m_0,\theta_0>0$ and $\int_{\Td}m_0\dd=1$.
  For each fixed $T>0$, there exists a sufficiently small $\eta_0>0$ such that, if
  \[
    \|m_0-\bar m\|_{L^\infty}
    +\|\theta_0-\bar\theta\|_{L^\infty}
    +\|u_T-\bar u\|_{W^{1,\infty}}<\eta_0
  \]
  then (KMR) admits a classical solution satisfying \hypref{hyp:regularity} and
  \hypref{hyp:smallness}, and this solution satisfies the estimates of Proposition~\ref{prop:L2-turnpike_main}.
\end{proposition}

Finally, we discuss the meaning of \eqref{hyp:spectral-gap}.
Fix the harvester density at $m=\bar m$ and linearize the resource equation around $\bar\theta$.
This gives
\[
\partial_t\zeta-\mu\Delta\zeta
=
\left\{K-2\bar\theta-\varepsilon(\rho*\bar m)\right\}\zeta.
\]
Set $\psi:=\zeta/\bar\theta$. Under \eqref{hyp:spectral-gap},
the stationary resource equation and integration by parts give
\[
\frac12\frac{d}{dt}
\int_{\mathbb T^d}\bar\theta\psi^2\,dx
=
-Q_\theta(\psi)
\le
-\lambda_\theta
\int_{\mathbb T^d}\bar\theta\psi^2\,dx.
\]
Thus, $Q_\theta$ represents the dissipation of the weighted $L^2$ energy for this linearized equation.
Condition \eqref{hyp:spectral-gap} ensures exponential decay of $\|\sqrt{\bar\theta}\,\psi(t)\|_{L^2(\mathbb T^d)}$
at a rate of at least $\lambda_\theta$.

The paper is organized as follows.
Section~2 establishes the modified energy identity and the local dissipation estimate (Proposition~\ref{prop:local-energy-di}).
Section~3 studies sufficient conditions for \hypref{hyp:spectral-gap}.
Section~4 derives the weighted $L^2$ turnpike estimates (Proposition~\ref{prop:L2-turnpike_main}).
Section~5 treats the $L^\infty$ turnpike estimate and the improvement of its constant.
It then uses these estimates to prove the existence of solutions near a stationary solution (Proposition~\ref{prop_exsist}).

\section{Dissipation estimate for the modified energy}
{
In this section, we compute
the time derivative of the modified energy $E$ introduced in Section~1
and prove the dissipation estimate in Proposition~\ref{prop:local-energy-di}.
We first use the zero-flux identity from the stationary Fokker--Planck equation (Lemma~\ref{lem:fluz})
to derive an identity for the time derivative of the modified energy (Proposition~\ref{prop:E-equation}).
We then split this identity into the leading quadratic part and higher-order error terms.
Under assumptions \eqref{hyp:spectral-gap} and \eqref{hyp:smallness}, we absorb the error terms into the leading part
to obtain the dissipation estimate.
In addition,
Section~2.2 shows that the same argument applies to a more general reaction term in place of $\theta(K-\theta)$ in the resource equation.
}
\subsection{Dissipation estimate for the modified energy}\label{subsec:proof-local-energy}
Let $(u,m,\theta)$ be a solution of (KMR), and let
$(\bar u,\bar m,\bar\theta,\lambda)$ be a solution of (sKMR).
We define the differences between these solutions by
\[
  v(t,x)=u(t,x)-\bar u(x)-\lambda(T-t),
  \qquad
  n(t,x)=m(t,x)-\bar m(x),
  \qquad
  z(t,x)=\theta(t,x)-\bar\theta(x).
\]
For convenience, we also set
$
  r:=\frac{n}{\bar m},
  \quad
  q:=\frac{z}{\bar\theta},
 \quad
  A:=\nabla v-\nu\nabla r.
$

We define the relative entropy of the resource density $\theta$ by
\[
  H_\theta(t)
  :=
  \int_{\Td}
  \bar\theta\{q-\log(1+q)\}\dd.
\]

The definition gives $H_{\theta}\ge0$ whenever $\theta, \bar{\theta}>0$.
Using this relative entropy, we define the modified energy $E$ by
\[
  E(t)
  :=
  -\int_{\Td}n v\dd
  +\nu\int_{\Td}\frac{n^2}{\bar m}\dd
  +\frac{1}{\eps}H_\theta(t).
\]

We record the following preliminary observations.
{
We first prove the zero-flux identity from the stationary Fokker--Planck equation.}
\begin{lemma}[Zero-flux identity]\label{lem:fluz}
\label{lem:perturbation-equations}
Assume (\ref{hyp:regularity}).
Then
\begin{equation}
 \nu\nabla\bar m=\bar m\nabla\bar u.
\end{equation}

\end{lemma}
\begin{proof}
We can write \eqref{eq:sKMR-FP} as
$
  \operatorname{div}\{\bar m\nabla\bar u-\nu\nabla\bar m\}=0
$.

Setting $\psi:=\bar m e^{-\bar u/\nu}$ gives
\[
  \bar m\nabla\bar u-\nu\nabla\bar m
  =
  -\nu e^{\bar u/\nu}\nabla\psi.
\]
Hence
$
  \operatorname{div}\left(e^{\bar u/\nu}\nabla\psi\right)=0
$.
Multiplying this equality by $\psi$ and integrating over $\mathbb T^d$, with no boundary terms, yields
$
  0
  =
  -\int_{\mathbb T^d}
  e^{\bar u/\nu}|\nabla\psi|^2\,dx
$.
Thus $\nabla\psi=0$, which proves the claim.
\end{proof}

\begin{remark}\label{rem:mass}
Let $m, \bar{m}$ be classical solutions of \eqref{eq:KMR-FP} and \eqref{eq:sKMR-FP}, respectively.
The divergence theorem gives
$\int_{\mathbb T^d}m\,dx=\int_{\mathbb T^d} m_0\,dx$.
Thus, by (\ref{hyp:mass}) and the additional conditions for (sKMR), $n=m-\bar{m}$ satisfies
$
    \int_{\Td}n\dd=0.
$
Since $n=\bar{m}r$, we also have
$
    \int_{\Td}\bar{m}r\dd=0.
$

\end{remark}

\begin{remark}\label{rem:commutative_*}
    Assume (\ref{hyp:even-kernel}).
    Then
  \[
    \int_{\Td} n(\rho*z)\dd
    =
    \int_{\Td} (\rho*n)z\dd.
  \]

  Indeed,
  \begin{align*}
    \int_{\Td} n(x)(\rho*z)(x)\,dx
    &=
    \int_{\Td}n(x)
    \left\{
      \int_{\Td}\rho(x-y)z(y)\,dy
    \right\}
    dx                                   \\
    &=
    \int_{\Td}z(y)
    \left\{
      \int_{\Td}\rho(y-x)n(x)\,dx
    \right\}
    dy                                      \\
    &=
    \int_{\Td}(\rho*n)(y)z(y)\,dy .
  \end{align*}
  Here we used $\rho(x-y)=\rho(y-x)$ in the penultimate equality.
\end{remark}

{
We now compute the time derivative of the modified energy $E$.
}
\begin{proposition}[Chain rule for the modified energy]
\label{prop:E-equation}
Assume \eqref{hyp:even-kernel}, \eqref{hyp:regularity}, and $\theta>0$.
Then the following identity holds.
\begin{align}
  E'(t)
  &=-\int_{\Td}\bar m\abs{\nabla v-\nu\nabla r}^2\dd
    -\nu^2\int_{\Td}\bar m\abs{\nabla r}^2\dd                           \notag\\
  &\quad
    -\frac{\mu}{\eps}
    \int_{\Td}
    \bar\theta\frac{\abs{\nabla q}^2}{(1+q)^2}\dd
    +
    \frac{\mu}{\eps}
    \int_{\Td}
    \frac{q}{1+q}\nabla\bar\theta\cdot\nabla q\dd                        \notag\\
  &\quad
    -\frac1\eps
    \int_{\Td}\bar\theta^2q^2\dd
    -\frac12\int_{\Td}n\abs{\nabla v}^2\dd
    +2\nu\int_{\Td}n\nabla r\cdot\nabla v\dd .
  \label{eq:energy-identity}
\end{align}
\end{proposition}

\begin{proof}
We compute the time derivative of the modified energy. Subtracting (sKMR) from (KMR) and rearranging gives
  \begin{align}
    -\partial_t v-\nu\Delta v-\nabla\bar u\cdot\nabla v
    &= \rho*z+\frac12\abs{\nabla v}^2,\label{eq:p1}\\
    \partial_t n
    &= \operatorname{div}\{\bar m(\nu\nabla r-\nabla v)\}
       -\operatorname{div}(n\nabla v),\label{eq:p2}\\
    \partial_t\theta
    &=
    {\mu\bar\theta\Delta q
    +2\mu\nabla\bar\theta\cdot\nabla q
    -\bar\theta^2q(1+q)
    -\eps\bar\theta (1+q)(\rho*n)}\label{eq:p3}.
  \end{align}
Here we used Lemma~\ref{lem:fluz} in the second equation
and $\theta={(1+q)}\bar\theta$ in the third.

Using the equations \eqref{eq:p1} for $v$ and \eqref{eq:p2} for $n$, we obtain
\begin{align*}
  \int_{\Td}(\partial_t n)v\dd
  &=-\int_{\Td}\bar m(\nu\nabla r-\nabla v)\cdot\nabla v\dd
    +\int_{\Td}n\abs{\nabla v}^2\dd                                      \\
  &=-\nu\int_{\Td}\bar m\nabla r\cdot\nabla v\dd
    +\int_{\Td}\bar m\abs{\nabla v}^2\dd
    +\int_{\Td}n\abs{\nabla v}^2\dd,
\end{align*}
On the other hand,
\begin{align*}
  \int_{\Td} n\,\partial_t v\dd
  &=
  \int_{\Td}
  n\left(
  -\nu\Delta v-\nabla\bar u\cdot\nabla v
  -\rho*z-\frac12\abs{\nabla v}^2
  \right)\dd                                      \\
  &=
  \nu\int_{\Td}\bar m\nabla r\cdot\nabla v\dd
  -\int_{\Td}n(\rho*z)\dd
  -\frac12\int_{\Td}n\abs{\nabla v}^2\dd .
\end{align*}
Here we used
the following calculation based on Lemma~\ref{lem:fluz},
\begin{align*}
\int_{\Td}n\Delta v\dd
&= \int_{\Td}-\nabla n\cdot \nabla v \dd\\
&= \int_{\Td}-\nabla(\bar m r)\cdot \nabla v \dd\\
&= \int_{\Td}-\bar m \nabla r\cdot \nabla v
   - r\nabla \bar m\cdot \nabla v \dd\\
&= \int_{\Td}-\bar m \nabla r\cdot \nabla v
   - \frac{1}{\nu} r\bar m \nabla \bar u\cdot \nabla v \dd\\
&= \int_{\Td}-\bar m \nabla r\cdot \nabla v
   - \frac{1}{\nu} n\nabla \bar u\cdot \nabla v \dd.
\end{align*}

By Remark~\ref{rem:commutative_*},
$
  \int_{\Td}n(\rho*z)\dd
  =
  \int_{\Td}(\rho*n)z\dd
$
and hence
\[
  \frac{d}{dt}\int_{\Td}nv\dd
  =
  \int_{\Td}\bar m\abs{\nabla v}^2\dd
  +\frac12\int_{\Td}n\abs{\nabla v}^2\dd
  -\int_{\Td}(\rho*n)z\dd .
\]

Next, since $\int_{\Td}n^2/\bar m\dd=\int_{\Td}\bar m r^2\dd$, we have
\begin{align*}
  \frac{d}{dt}\int_{\Td}\frac{n^2}{\bar m}\dd
  &=2\int_{\Td}r\,\partial_t n\dd                                      \\
  &=2\int_{\Td}
    r\,\operatorname{div}\{\bar m(\nu\nabla r-\nabla v)\}\dd
    -2\int_{\Td}r\,\operatorname{div}(n\nabla v)\dd                     \\
  &=-2\nu\int_{\Td}\bar m\abs{\nabla r}^2\dd
    +2\int_{\Td}\bar m\nabla r\cdot\nabla v\dd
    +2\int_{\Td}n\nabla r\cdot\nabla v\dd .
\end{align*}
It follows that
\begin{align*}
  &\frac{d}{dt}\left(
    -\int_{\Td}nv\dd
    +\nu\int_{\Td}\frac{n^2}{\bar m}\dd
  \right)                                                   \\
  &=-\int_{\Td}\bar m\abs{\nabla v}^2\dd
    +2\nu\int_{\Td}\bar m\nabla r\cdot\nabla v\dd
    -2\nu^2\int_{\Td}\bar m\abs{\nabla r}^2\dd           \\
  &\quad
    +\int_{\Td}(\rho*n)z\dd
    -\frac12\int_{\Td}n\abs{\nabla v}^2\dd
    +2\nu\int_{\Td}n\nabla r\cdot\nabla v\dd              \\
  &=-\int_{\Td}\bar m\abs{\nabla v-\nu\nabla r}^2\dd
    -\nu^2\int_{\Td}\bar m\abs{\nabla r}^2\dd             \\
  &\quad
    +\int_{\Td}(\rho*n)z\dd
    -\frac12\int_{\Td}n\abs{\nabla v}^2\dd
    +2\nu\int_{\Td}n\nabla r\cdot\nabla v\dd .
\end{align*}

Next, we differentiate the relative entropy of \(\theta\).
By definition,
\[
H_\theta(t)=\int_{\mathbb T^d}\bar\theta {\{q-\log(1+q)\}}\,dx.
\]

Using $\partial_t\theta=\bar\theta\,\partial_t q$ and \eqref{eq:p3},
we obtain
\begin{align*}
\frac{d}{dt}H_\theta(t)
&=
\int_{\mathbb T^d}
\bar\theta {\frac{q}{1+q}}\partial_t q\,dx\\
  &=
  \int_{\Td}\left({\frac{q}{1+q}}\right)\partial_t\theta\dd                                      \\
  &=
  \mu\int_{\Td}
  \left({\frac{q}{1+q}}\right)
  \left(
    \bar\theta\Delta {q}
    +2\nabla\bar\theta\cdot\nabla {q}
  \right)\dd                        
  -\int_{\Td}\bar\theta^2{q}^2\dd
  -\eps\int_{\Td}\bar\theta q(\rho*n)\dd .
\end{align*}
Since $z=\bar\theta{q}$, the last term is
$
  -\eps\int_{\Td}\bar\theta{q}(\rho*n)\dd
  =
  -\eps\int_{\Td}z(\rho*n)\dd
$.

We integrate the diffusion term by parts. First,
\[
\begin{aligned}
&\int_{\mathbb T^d}
\left({\frac{q}{1+q}}\right)
\left(
\bar\theta\Delta {q}+2\nabla\bar\theta\cdot\nabla {q}
\right)\,dx  \\
&\quad =
\int_{\mathbb T^d}
\left({\frac{q}{1+q}}\right)\bar\theta\Delta {q}\,dx
+
2\int_{\mathbb T^d}
{\frac{q}{1+q}}
\nabla\bar\theta\cdot\nabla {q}\,dx .
\end{aligned}
\]
Integrating the first diffusion term by parts on \(\mathbb T^d\), with no boundary terms, gives
\[
\begin{aligned}
\int_{\mathbb T^d}
\left({\frac{q}{1+q}}\right)\bar\theta\Delta {q}\,dx
&=
-\int_{\mathbb T^d}
\nabla\left\{
\left({\frac{q}{1+q}}\right)\bar\theta
\right\}\cdot\nabla {q}\,dx  \\
&=
-\int_{\mathbb T^d}
\left\{
\bar\theta\nabla\left({\frac{q}{1+q}}\right)
+
{\frac{q}{1+q}}\nabla\bar\theta
\right\}\cdot\nabla {q}\,dx .
\end{aligned}
\]
Since
$
\nabla\left({\frac{q}{1+q}}\right)
=
\frac{\nabla {q}}{{(1+q)^2}}
$,
we obtain
\[
\begin{aligned}
\frac{d}{dt}H_\theta(t)
&=
-\mu\int_{\mathbb T^d}
\bar\theta
\frac{|\nabla q|^2}{(1+q)^2}\,dx
+
\mu\int_{\mathbb T^d}
\frac{q}{1+q}
\nabla\bar\theta\cdot\nabla q\,dx  \\
&\quad
-
\int_{\mathbb T^d}
\bar\theta^2q^2\,dx
-
\varepsilon
\int_{\mathbb T^d}
z(\rho*n)\,dx .
\end{aligned}
\]
Combining these identities gives
\begin{align}
  E'(t)
  &=-\int_{\Td}\bar m\abs{\nabla v-\nu\nabla r}^2\dd
    -\nu^2\int_{\Td}\bar m\abs{\nabla r}^2\dd                           \notag\\
  &\quad
    -\frac{\mu}{\eps}
    \int_{\Td}
    \bar\theta\frac{\abs{\nabla q}^2}{(1+q)^2}\dd
    +
    \frac{\mu}{\eps}
    \int_{\Td}
    \frac{q}{1+q}\nabla\bar\theta\cdot\nabla q\dd                        \notag\\
  &\quad
    -\frac1\eps
    \int_{\Td}\bar\theta^2q^2\dd
    -\frac12\int_{\Td}n\abs{\nabla v}^2\dd
    +2\nu\int_{\Td}n\nabla r\cdot\nabla v\dd .
\end{align}
This is the required energy identity.
\end{proof}

{
We now split
the right-hand side of this identity
into the leading quadratic part and error terms of order at least three in the perturbations.
We introduce the following quantities.
}
The flux dissipation associated with \eqref{eq:p2} is
\[
  D_{\mathrm{flux}}(t)
  =
  \int_{\Td}\bar m\abs{A}^2\dd
  =
  \int_{\Td}\bar m\abs{\nabla v-\nu\nabla r}^2\dd.
\]

The resource quadratic form is
\[
  Q_\theta(\varphi)
  =
  \mu\int_{\Td}\bar\theta\abs{\nabla\varphi}^2\dd
  -\mu\int_{\Td}\varphi\nabla\bar\theta\cdot\nabla\varphi\dd
  +\int_{\Td}\bar\theta^2\varphi^2\dd.
\]
Combining this form with the density-gradient term gives the quadratic form
\[
  \mathcal{Q}(r,q)
  =
  \nu^2\int_{\Td}\bar m\abs{\nabla r}^2\dd
  +\frac{1}{\eps}Q_\theta(q).
\]

With this notation,
the derivative of the modified energy takes the form
\begin{equation}
  E'(t)
  =
  -D_{\mathrm{flux}}(t)
  -\mathcal{Q}(r,q)
  +R_m(t)+R_\theta(t)
  \label{eq:energy-decomposition}
\end{equation}
where the remaining error terms are
\[
  R_m(t)
  =
  -\frac12\int_{\Td}n\abs{\nabla v}^2\dd
  +2\nu\int_{\Td}n\nabla r\cdot\nabla v\dd,
\]
\begin{align*}
  R_\theta(t)
  ={}&
  \frac{\mu}{\eps}
  \int_{\Td}\bar\theta
  \left\{
    1-\frac{1}{(1+q)^2}
  \right\}
  \abs{\nabla q}^2\dd                                      
  +
  \frac{\mu}{\eps}
  \int_{\Td}
  \left\{
    \frac{q}{1+q}-q
  \right\}
  \nabla\bar\theta\cdot\nabla q\dd.
\end{align*}

Below we write $R_{\mathrm{total}}:=R_m+R_\theta$.

{
We now
prove Proposition~\ref{prop:local-energy-di}.
We first
estimate the error terms $R_m$ and $R_\theta$ separately
and absorb them into the leading terms $D_{\mathrm{flux}}+\mathcal{Q}$ under \eqref{hyp:regularity}--\eqref{hyp:smallness}.
This gives the dissipation estimate \eqref{eq:main-dissipation} for the modified energy.
We then
bound the absolute value of the energy by the leading part $D_{\mathrm{flux}}+\mathcal{Q}$
to obtain \eqref{eq:main-energy-control}.
}
\begin{proof}[Proof of Proposition~\ref{prop:local-energy-di}]
We first bound $R_m$ by the leading terms $D_{\mathrm{flux}}+\mathcal{Q}$.
Using $n=\bar m r$ and $\nabla v=A+\nu\nabla r$, we obtain
\begin{align*}
  R_m
  &=
  -\frac12\int_{\Td}\bar m r\abs{A+\nu\nabla r}^2\dd
  +2\nu\int_{\Td}\bar m r\nabla r\cdot(A+\nu\nabla r)\dd                  \\
  &=
  -\frac12\int_{\Td}\bar m r\abs{A}^2\dd
  +\nu\int_{\Td}\bar m r A\cdot\nabla r\dd
  +\frac{3\nu^2}{2}\int_{\Td}\bar m r\abs{\nabla r}^2\dd .
\end{align*}
Young's inequality then gives
\begin{align}
  \abs{R_m}
  &\leq
  2\norm[L^\infty]{r}
  \left\{
    \int_{\Td}\bar m\abs{A}^2\dd
    +\nu^2\int_{\Td}\bar m\abs{\nabla r}^2\dd
  \right\}                                                   \notag\\
  &\leq
  2\norm[L^\infty]{r}
  \left\{
    D_{\mathrm{flux}}+\mathcal{Q}(r,q)
  \right\} .
  \label{eq:Rm-estimate}
\end{align}
The last inequality uses $Q_{\theta}\ge0$, which follows from \eqref{hyp:spectral-gap}.

Next, we bound $R_\theta$ by $D_{\mathrm{flux}}+\mathcal{Q}$.
Since $\norm[L^\infty]{q}\leq\delta\leq1/2$, we have $1+q\geq1/2$ for $x$. 
Thus
\[
  \left\lvert
    1-\frac{1}{(1+q)^2}
  \right\rvert
  \leq
  6\abs{q}.
\]
Also,
$
  \frac{q}{1+q}-q
  =
  -\frac{q^2}{1+q}
$
gives
\[
  \left\lvert
    \frac{q}{1+q}-q
  \right\rvert
  \leq
  2\abs{q}^2
  \leq
  2\norm[L^\infty]{q}\abs{q}.
\]

Hence
the first term of \(R_\theta\) satisfies
\[
\begin{aligned}
\left|
\frac{\mu}{\eps}
\int_{\Td}
\bar\theta
\left\{
1-\frac{1}{(1+q)^2}
\right\}
\abs{\nabla q}^2\dd
\right|
&\leq
6\norm[L^\infty]{q}
\frac{\mu}{\eps}
\int_{\Td}\bar\theta\abs{\nabla q}^2\dd.
\end{aligned}
\]
For the second term, we have
\[
\begin{aligned}
\left|
\frac{\mu}{\eps}
\int_{\Td}
\left\{
\frac{q}{1+q}-q
\right\}
\nabla\bar\theta\cdot\nabla q\dd
\right|
&\leq
2\norm[L^\infty]{q}
\frac{\mu}{\eps}
\int_{\Td}
\abs{q}\abs{\nabla\bar\theta}\abs{\nabla q}\dd .
\end{aligned}
\]

The positivity of \(\bar\theta\) and boundedness of \(\nabla\bar\theta\) in (\ref{hyp:regularity}) give
$
\frac{|\nabla \bar{\theta}|^2}{\bar{\theta}}\le \|\frac{\nabla \bar \theta}{\bar \theta}\|_{L^\infty} ^2 \bar{\theta}
$.

Thus, Young's inequality yields
\[
\begin{aligned}
\mu\abs{q}\abs{\nabla\bar\theta}\abs{\nabla q}
&=
\left(\sqrt{\mu\bar\theta}\abs{\nabla q}\right)
\left(
\sqrt{\mu}\frac{\abs{\nabla\bar\theta}}{\sqrt{\bar\theta}}\abs{q}
\right)  \\
&\leq
\frac12\mu\bar\theta\abs{\nabla q}^2
+
\frac{\mu}{2} \left\|\frac{\nabla \bar \theta}{\bar \theta}\right\|_{L^\infty} ^2 \bar\theta q^2 .
\end{aligned}
\]
Consequently,
\[
\begin{aligned}
\left|
\frac{\mu}{\eps}
\int_{\Td}
\left\{
\frac{q}{1+q}-q
\right\}
\nabla\bar\theta\cdot\nabla q\dd
\right|
&\leq
\norm[L^\infty]{q}
\frac{1}{\eps}
\left\{
\mu\int_{\Td}\bar\theta\abs{\nabla q}^2\dd
+
\mu \left\|\frac{\nabla \bar \theta}{\bar \theta}\right\|_{L^\infty} ^2\int_{\Td}\bar\theta q^2\dd
\right\}.
\end{aligned}
\]
Combining these two estimates, we obtain
\[
  \abs{R_\theta}
  \leq
  \norm[L^\infty]{q}
  \frac{1}{\eps}
  \left\{
    7 \mu\int_{\Td}\bar\theta\abs{\nabla q}^2\dd
    +\mu \left\|\frac{\nabla \bar \theta}{\bar \theta}\right\|_{L^\infty} ^2 \int_{\Td}\bar\theta q^2\dd
  \right\}.
\]

We bound the integrals on the right by $D_{\mathrm{flux}}+\mathcal{Q}$.
Integration by parts gives, for every $\varphi\in H^1(\mathbb T^d)$,
\[
-\mu\int_{\mathbb T^d}\varphi\nabla\bar\theta\cdot\nabla\varphi\,dx
=
-\frac{\mu}{2}\int_{\mathbb T^d}\nabla\bar\theta\cdot\nabla(\varphi^2)\,dx
=
\frac{\mu}{2}\int_{\mathbb T^d}(\Delta\bar\theta)\varphi^2\,dx.
\]
Hence
\[
Q_\theta(\varphi)
=
\mu\int_{\mathbb T^d}\bar\theta|\nabla\varphi|^2\,dx
+
\int_{\mathbb T^d}V_\theta\varphi^2\,dx,
\qquad
V_\theta
:=
\bar\theta^2+\frac{\mu}{2}\Delta\bar\theta .
\]
By \eqref{eq:sKMR-theta},
$
\mu\Delta\bar\theta
=
-\bar\theta K+\bar\theta^2
+\varepsilon\bar\theta(\rho*\bar m)
$
and therefore
\[
V_\theta
=
\bar\theta^2
+\frac12\{-\bar\theta K+\bar\theta^2
+\varepsilon\bar\theta(\rho*\bar m)\}
=
\bar\theta
\left(
\frac32\bar\theta
-\frac12 K
+\frac{\varepsilon}{2}(\rho*\bar m)
\right).
\]
The boundedness of $\bar\theta$, $K$, and $\rho*\bar m$ gives a constant $C_V>0$ such that
\[
|V_\theta(x)|
\le
C_V\bar\theta(x)
\qquad
(x\in\mathbb T^d).
\]
For example, Young's inequality and Remark~\ref{rem:mass} allow us to take
\[
C_V
=
\frac32\|\bar\theta\|_{L^\infty}
+
\frac12\|K\|_{L^\infty}
+
\frac{\varepsilon}{2}\|\rho\|_{L^\infty}.
\]

It follows that
\[
Q_\theta(\varphi)
=
\mu\int_{\mathbb T^d}\bar\theta|\nabla\varphi|^2\,dx
+
\int_{\mathbb T^d}V_\theta\varphi^2\,dx
\ge
\mu\int_{\mathbb T^d}\bar\theta|\nabla\varphi|^2\,dx
-
C_V\int_{\mathbb T^d}\bar\theta\varphi^2\,dx .
\]
Thus
\[
\mu\int_{\mathbb T^d}\bar\theta|\nabla\varphi|^2\,dx
\le
Q_\theta(\varphi)
+
C_V\int_{\mathbb T^d}\bar\theta\varphi^2\,dx .
\]
Moreover, assumption \eqref{hyp:spectral-gap} gives
$
\int_{\mathbb T^d}\bar\theta\varphi^2\,dx
\le
\lambda_\theta^{-1}Q_\theta(\varphi)
$
and hence
\[
\mu\int_{\mathbb T^d}\bar\theta|\nabla\varphi|^2\,dx
\le
\left(1+\frac{C_V}{\lambda_\theta}\right)Q_\theta(\varphi).
\]
Using this estimate and \eqref{hyp:spectral-gap} once more, we obtain
\begin{equation}
7\mu\int_{\mathbb T^d}\bar\theta|\nabla\varphi|^2\,dx
+
\mu\left\|\frac{\nabla \bar \theta}{\bar \theta}\right\|_{L^\infty} ^2\int_{\mathbb T^d}\bar\theta\varphi^2\,dx
\le
\left(\frac{\mu}{\lambda_\theta}\left\|\frac{\nabla \bar \theta}{\bar \theta}\right\|_{L^\infty} ^2+7(1+\frac{C_V}{\lambda_{\theta}})\right) Q_\theta(\varphi)
\qquad
(\varphi\in H^1(\mathbb T^d)).
\end{equation}

We can now estimate
\begin{align}
  \abs{R_\theta}
  &\leq
  \norm[L^\infty]{q}
  \frac{1}{\eps}
  \left\{
    7\mu\int_{\Td}\bar\theta\abs{\nabla q}^2\dd
    +\mu \left\|\frac{\nabla \bar \theta}{\bar \theta}\right\|_{L^\infty} ^2\int_{\Td}\bar\theta q^2\dd
  \right\}                                                   \notag\\
  &\leq
  \left(\frac{\mu}{\lambda_\theta}\left\|\frac{\nabla \bar \theta}{\bar \theta}\right\|_{L^\infty} ^2+7(1+\frac{C_V}{\lambda_{\theta}})\right)
  \norm[L^\infty]{q}
  \frac{1}{\eps}Q_\theta(q)                                  \notag\\
  &\leq
  \left(\frac{\mu}{\lambda_\theta}\left\|\frac{\nabla \bar \theta}{\bar \theta}\right\|_{L^\infty} ^2+7(1+\frac{C_V}{\lambda_{\theta}})\right)
  \norm[L^\infty]{q}\mathcal{Q}(r,q) .
  \label{eq:Rtheta-estimate}
\end{align}

Since $D_{\mathrm{flux}}$ is nonnegative, \eqref{eq:Rm-estimate} and \eqref{eq:Rtheta-estimate} yield
\[
  \abs{R_{\mathrm{total}}}
  \leq
  \max\left\{2, \left(\frac{\mu}{\lambda_\theta}\left\|\frac{\nabla \bar \theta}{\bar \theta}\right\|_{L^\infty} ^2+7(1+\frac{C_V}{\lambda_{\theta}})\right) \right\}
  \left(
    \norm[L^\infty]{r}+
    \norm[L^\infty]{q}
  \right)
  \left\{
    D_{\mathrm{flux}}+\mathcal{Q}(r,q)
  \right\}.
\]
By \eqref{hyp:smallness}, this implies
\[
  \abs{R_{\mathrm{total}}}
  \leq
  \frac12
  \left\{
    D_{\mathrm{flux}}+\mathcal{Q}(r,q)
  \right\}.
\]
Substituting this estimate into \eqref{eq:energy-decomposition}, we obtain
\[
  -E'(t)
  \geq
  \frac12
  \left\{
    D_{\mathrm{flux}}(t)+\mathcal{Q}(r(t),q(t))
  \right\}.
\]
This proves \eqref{eq:main-dissipation}.

Finally, we bound $\abs{E}$ itself by the leading terms.
For the second term in $E$, we have
$
  \nu\int_{\Td}\frac{n^2}{\bar m}\dd
  =
  \nu\int_{\Td}\bar m r^2\dd
$.

By Remark~\ref{rem:mass} and the weighted Poincar\'e--Wirtinger inequality,
\[
  \int_{\Td}\bar m r^2\dd
  \leq
  C_P\int_{\Td}\bar m\abs{\nabla r}^2\dd .
\]
Since (\ref{hyp:spectral-gap}) implies that $Q_{\theta}$ is nonnegative, we obtain
\begin{equation}\label{eq:est_Q_by_P}
  \int_{\Td}\bar m r^2\dd
  \leq
  C \mathcal{Q}(r,q).
\end{equation}

Next, recall that $\norm[L^\infty]{q}\leq1/2$ and that, for $|x|\le \frac{1}{2}$, there is a constant $C>0$ such that $0 \le x-\log(1+x)\le C x^2$. Hence the third term in $E$ satisfies
\[
  H_\theta(t)
  \leq
  C\int_{\Td}\bar\theta q^2\dd
  \leq
  C Q_\theta(q).
\]
It follows that
\[
  \frac{1}{\eps}H_\theta(t)
  \leq
  C \mathcal{Q}(r,q) .
\]

We finally estimate the first term in $E$. Set
$\baravg{v}
  :=
  \frac{\int_{\Td}\bar m v\dd}{\int_{\Td}\bar m\dd}
$.
By Remark~\ref{rem:mass} and the weighted Poincar\'e--Wirtinger inequality,
\[
  \int_{\Td}\bar m\abs{v-\baravg{v}}^2\dd
  \leq
  C_P\int_{\Td}\bar m\abs{\nabla v}^2\dd .
\]
The Cauchy--Schwarz inequality then gives
\[
  \left\lvert\int_{\Td}nv\dd\right\rvert
  =
  \left\lvert
  \int_{\Td} \bar m r\left(v-\baravg{v}\right)\dd
  \right\rvert
  \leq
  \sqrt{C_P}
  \left(\int_{\Td}\bar m r^2\dd\right)^{1/2}
  \left(\int_{\Td}\bar m\abs{\nabla v}^2\dd\right)^{1/2}.
\]
Using $\nabla v=A+\nu\nabla r$ and the nonnegativity of $Q_{\theta}$ from \eqref{hyp:spectral-gap}, we obtain
\begin{equation*}
  \int_{\Td}\bar m\abs{\nabla v}^2\dd
  \leq
  C
  \left\{
    D_{\mathrm{flux}}+\nu^2\int_{\Td}\bar m\abs{\nabla r}^2\dd
  \right\}
  \leq
  C
  \left\{
    D_{\mathrm{flux}}+\mathcal{Q}(r,q)
  \right\} .
\end{equation*}
Combining this gradient estimate with \eqref{eq:est_Q_by_P} gives
\[
  \left\lvert\int_{\Td}nv\dd\right\rvert
  \leq
  C
  \mathcal{Q}(r,q)^{1/2}
  \left\{
    D_{\mathrm{flux}}+\mathcal{Q}(r,q)
  \right\}^{1/2}
  \leq
  C
  \left\{
    D_{\mathrm{flux}}+\mathcal{Q}(r,q)
  \right\} .
\]

Combining the three estimates, we find a positive constant $C_E$ such that
\[
  \abs{E(t)}
  \leq
  C_E
  \left\{
    D_{\mathrm{flux}}(t)+\mathcal{Q}(r(t),q(t))
  \right\}.
\]
This proves \eqref{eq:main-energy-control}. Finally, \eqref{eq:main-dissipation} and
\eqref{eq:main-energy-control} give
\[
  -E'(t)
  \geq
  \frac{1}{2C_E}\abs{E(t)}.
\]
Setting $\sigma=\frac{1}{2C_E}$ yields \eqref{eq:main-signed-energy}.
\end{proof}

\begin{remark}[Integral form]
  For every $0<t_1<t_2<T$, integration of \eqref{eq:main-dissipation} gives
  \[
    \frac{1}{2}\int_{t_1}^{t_2}
    \left\{
      D_{\mathrm{flux}}(s)+\mathcal{Q}(r(s),q(s))
    \right\}\,ds
    \leq
    E(t_1)-E(t_2).
  \]
  This integral form can be combined with information from the initial and terminal times to derive turnpike estimates.
\end{remark}
\begingroup
\subsection{Extension to general reaction terms}\label{subsec:general-reaction}
In this subsection, we replace the reaction term $\theta(K-\theta)$ in the resource equation by a general reaction term
$G:\Td\times(0,\infty)\to\mathbb R$
and extend Proposition~\ref{prop:local-energy-di}.
The first variable $x\in\Td$ in $G(x,s)$ denotes position, and the second variable $s>0$ denotes the resource density.

We retain \eqref{eq:KMR-HJ}, \eqref{eq:KMR-FP}, their corresponding stationary equations,
the initial and terminal conditions, and the normalization of the stationary solution.
We replace the resource equation and its stationary counterpart by
\begin{align}
  \partial_t\theta-\mu\Delta\theta
  &=G(x,\theta)-\varepsilon(\rho*m)\theta,
  \tag{\textup{KMR-3}$_G$}\label{eq:G-resource}\\
  -\mu\Delta\bar\theta
  &=G(x,\bar\theta)-\varepsilon(\rho*\bar m)\bar\theta
  \tag{\textup{sKMR-3}$_G$}\label{eq:G-stationary-resource}.
\end{align}
The perturbations $v,n,z,r,q,A$ and the modified energy $E$
are defined by the same formulas for solutions of the modified systems.
The relative entropy and the coupling structure of the harvesting term are preserved.
We therefore modify the quadratic form containing the quadratic part of the reaction term and the conditions used to absorb the remainder terms.

We replace \hypref{hyp:regularity}, \hypref{hyp:spectral-gap}, and \hypref{hyp:smallness}
by the following conditions.

\medskip\noindent
\textbf{(H3$_G$) Classical solutions and regularity of the reaction term.}
The functions $(u,m,\theta)$ and $(\bar u,\bar m,\bar\theta,\lambda)$ are classical solutions
of the time-dependent system \eqref{eq:KMR-HJ}, \eqref{eq:KMR-FP}, \eqref{eq:G-resource}
and the stationary system \eqref{eq:sKMR-HJ}, \eqref{eq:sKMR-FP}, \eqref{eq:G-stationary-resource}, respectively.
For constants $C_m,C_\theta\geq1$, they satisfy
\[
  C_m^{-1}\leq\bar m(x)\leq C_m,\qquad
  C_\theta^{-1}\leq\bar\theta(x)\leq C_\theta,
  \qquad x\in\Td.
\]
We also assume the regularity condition
\[
  G\in C^2(\Td\times(0,\infty)).
\]
For example, this condition holds for the logistic reaction $G(x,s)=s(K(x)-s)$ with $K\in C^2(\Td)$,
and for the Nagumo cubic $G(x,s)=s(1-s)(s-a)$ with $0<a<1$.

For $s>0$, set $g(x,s):=G(x,s)/s$, and define
\[
  a_G(x):=G(x,\bar\theta(x))
    -\bar\theta(x)\partial_sG(x,\bar\theta(x))
    =-\bar\theta(x)^2\partial_sg(x,\bar\theta(x))
\]
and
\[
  Q_G(\varphi):=
  \mu\int_{\Td}\bar\theta|\nabla\varphi|^2\dd
  -\mu\int_{\Td}\varphi\nabla\bar\theta\cdot\nabla\varphi\dd
  +\int_{\Td}a_G\varphi^2\dd.
\]
The corresponding quadratic form is
\[
  \mathcal{Q}_G(r,q):=
  \nu^2\int_{\Td}\bar m|\nabla r|^2\dd
  +\frac1\varepsilon Q_G(q).
\]
For $G(x,s)=s(K(x)-s)$, we have $a_G=\bar\theta^2$,
so $Q_G$ and $\mathcal{Q}_G$ coincide with $Q_\theta$ and $\mathcal{Q}$, respectively.

\medskip\noindent
\textbf{(H5$_G$) Coercivity of $Q_G$.}
There exists $\lambda_G>0$ such that, for every $\varphi\in H^1(\Td)$,
\[
  Q_G(\varphi)\geq\lambda_G\int_{\Td}\bar\theta\varphi^2\dd
\]

\medskip\noindent
\textbf{(H6$_G$) Smallness of the relative perturbations.}
In this subsection, we redefine the constants $C_V,C_R$ in \hypref{hyp:smallness} as follows.
First, set
\[
  M:=\sup_{\substack{x\in\Td\\
    \bar\theta(x)/2\leq s\leq3\bar\theta(x)/2}}
    |\partial_s^2g(x,s)|
\]
and define
\[
  C_V:=\left\|\frac{G(\cdot,\bar\theta)}{2\bar\theta}
    -\partial_sG(\cdot,\bar\theta)\right\|_{L^\infty}
    +\frac\varepsilon2\|\rho\|_{L^\infty},
\]
\[
  C_R:=\max\left\{2,\;
    7\left(1+\frac{C_V}{\lambda_G}\right)
    +\frac\mu{\lambda_G}
      \left\|\frac{\nabla\bar\theta}{\bar\theta}\right\|_{L^\infty}^2
    +\frac{M\|\bar\theta\|_{L^\infty}^2}{2\lambda_G}\right\}.
\]
There exists $\delta\in(0,1/2]$ such that
\[
  C_R\delta\leq\frac12,\qquad
  \|r(t)\|_{L^\infty}+\|q(t)\|_{L^\infty}\leq\delta,
  \qquad t\in[0,T].
\]

\begin{proposition}[Dissipation estimate for the modified energy with a general reaction term]
\label{prop:general-reaction}
Assume \hypref{hyp:even-kernel}, \hypref{hyp:mass}, and
(H3$_G$), (H5$_G$), (H6$_G$).
Then there exists a constant $C_{E,G}>0$, independent of $T$, such that,
for every $t\in[0,T]$,
\begin{align}
  -E'(t)&\geq\frac12
    \{D_{\mathrm{flux}}(t)+\mathcal{Q}_G(r(t),q(t))\},
    \label{eq:G-dissipation}\\
  |E(t)|&\leq C_{E,G}
    \{D_{\mathrm{flux}}(t)+\mathcal{Q}_G(r(t),q(t))\}.
    \label{eq:G-energy-control}
\end{align}
Consequently, setting $\sigma_G:=1/(2C_{E,G})>0$ gives
\begin{equation}
  E'(t)\leq-\sigma_G|E(t)|.
  \label{eq:G-signed-energy}
\end{equation}

\end{proposition}

\begin{proof}
We first compute the time derivative of the modified energy and bound the remainder terms by the leading part.
We then estimate the energy itself in terms of the leading part.

The Hamilton--Jacobi and Fokker--Planck equations are the same as in the original system.
Thus the calculation in the introduction gives
\[
\begin{aligned}
  E_2'(t)
  &=-D_{\mathrm{flux}}(t)-\nu^2\int_{\Td}\bar m|\nabla r|^2\dd
    +\int_{\Td}(\rho*n)z\dd+R_m(t),\\
  R_m(t)&=-\frac12\int_{\Td}n|\nabla v|^2\dd
    +2\nu\int_{\Td}n\nabla r\cdot\nabla v\dd.
\end{aligned}
\]
On the other hand, using $\theta=\bar\theta (1+q)$
and the stationary resource equation \eqref{eq:G-stationary-resource},
we can write \eqref{eq:G-resource} as
\[
\begin{aligned}
  \partial_t\theta
  ={}&\mu\bar\theta\Delta q+2\mu\nabla\bar\theta\cdot\nabla q+\bar\theta (1+q)\{g(x,\bar\theta (1+q))-g(x,\bar\theta)\}
    -\varepsilon\bar\theta (1+q)(\rho*n).
\end{aligned}
\]
Substituting this into
$H_\theta'=\int_{\Td}\frac{q}{1+q}\partial_t\theta\dd$ and integrating by parts yields
\[
\begin{aligned}
  H_\theta'(t)
  ={}&-\mu\int_{\Td}\bar\theta
        \frac{|\nabla q|^2}{(1+q)^2}\dd
      +\mu\int_{\Td}\frac q{1+q}\nabla\bar\theta\cdot\nabla q\dd\\
  &+\int_{\Td}z\{g(x,\bar\theta+z)-g(x,\bar\theta)\}\dd
      -\varepsilon\int_{\Td}z(\rho*n)\dd.
\end{aligned}
\]
Thus, in $E'=E_2'+\varepsilon^{-1}H_\theta'$,
the coupling terms cancel as in the logistic case.

To split the reaction contribution into its quadratic part and a remainder, set
\[
  R_G(t):=\frac1\varepsilon\int_{\Td}z
    \{g(x,\bar\theta+z)-g(x,\bar\theta)
      -\partial_sg(x,\bar\theta)z\}\dd.
\]
Since $z=\bar\theta q$, the definition of $a_G$ gives
\[
  \int_{\Td}z\{g(x,\bar\theta+z)-g(x,\bar\theta)\}\dd
  =-\int_{\Td}a_Gq^2\dd+\varepsilon R_G(t).
\]
The remainder from the diffusion terms is the same as before,
\[
\begin{aligned}
  R_\theta(t):={}&\frac\mu\varepsilon\int_{\Td}\bar\theta
    \left\{1-\frac1{(1+q)^2}\right\}|\nabla q|^2\dd+\frac\mu\varepsilon\int_{\Td}
      \left\{\frac q{1+q}-q\right\}\nabla\bar\theta\cdot\nabla q\dd.
\end{aligned}
\]
Combining these expressions, we obtain the energy identity
\begin{equation}
  E'=-D_{\mathrm{flux}}-\mathcal{Q}_G+R_m+R_\theta+R_G.
  \label{eq:G-energy-decomposition}
\end{equation}

We first estimate $R_m$. Using $n=\bar m r$ and $\nabla v=A+\nu\nabla r$, we have
\[
  R_m=-\frac12\int_{\Td}\bar m r|A|^2\dd
    +\nu\int_{\Td}\bar m rA\cdot\nabla r\dd
    +\frac{3\nu^2}{2}\int_{\Td}\bar m r|\nabla r|^2\dd.
\]
Young's inequality and $Q_G\geq0$, which follows from (H5$_G$), then give
\begin{equation}
  |R_m|\leq2\|r\|_{L^\infty}
    \{D_{\mathrm{flux}}+\mathcal{Q}_G\}.
  \label{eq:G-Rm-estimate}
\end{equation}

Next, we estimate $R_\theta$. Since $\|q\|_{L^\infty}\leq1/2$,
\[
  \left|1-\frac1{(1+q)^2}\right|\leq6|q|,\qquad
  \left|\frac q{1+q}-q\right|\leq2|q|^2.
\]
Using these inequalities and Young's inequality as before, we obtain
\[
  |R_\theta|\leq\frac{\|q\|_{L^\infty}}\varepsilon
    \left\{7\mu\int_{\Td}\bar\theta|\nabla q|^2\dd
      +\mu\left\|\frac{\nabla\bar\theta}{\bar\theta}\right\|_{L^\infty}^2\int_{\Td}\bar\theta q^2\dd\right\}.
\]
To bound the right-hand side by $Q_G$, we integrate by parts and write
\[
Q_G(\varphi)=\mu\int_{\Td}\bar\theta|\nabla\varphi|^2\dd
    +\int_{\Td}V_G\varphi^2\dd,\qquad
  V_G:=a_G+\frac\mu2\Delta\bar\theta.
\]
The stationary resource equation \eqref{eq:G-stationary-resource} gives
$
  \mu\Delta\bar\theta=-G(x,\bar\theta)
    +\varepsilon\bar\theta(\rho*\bar m).
$
Combining this with the definition of $a_G$ yields
\[
\begin{aligned}
  V_G
  &=\frac12G(x,\bar\theta)-\bar\theta\partial_sG(x,\bar\theta)
    +\frac\varepsilon2\bar\theta(\rho*\bar m)\\
  &=\bar\theta\left\{\frac{G(x,\bar\theta)}{2\bar\theta}
    -\partial_sG(x,\bar\theta)+\frac\varepsilon2(\rho*\bar m)\right\}.
\end{aligned}
\]
Since $\bar m\geq0$ and $\int_{\Td}\bar m\dd=1$,
we have $\|\rho*\bar m\|_{L^\infty}\leq\|\rho\|_{L^\infty}$, and hence
\[
\begin{aligned}
  |V_G(x)|
  &\leq\bar\theta(x)\left\{
    \left|\frac{G(x,\bar\theta(x))}{2\bar\theta(x)}
      -\partial_sG(x,\bar\theta(x))\right|
    +\frac\varepsilon2\|\rho\|_{L^\infty}\right\}
  \leq C_V\bar\theta(x).
\end{aligned}
\]
This estimate and (H5$_G$) imply
\[
\begin{aligned}
  \mu\int_{\Td}\bar\theta|\nabla\varphi|^2\dd
  &\leq Q_G(\varphi)+C_V\int_{\Td}\bar\theta\varphi^2\dd\leq\left(1+\frac{C_V}{\lambda_G}\right)Q_G(\varphi).
\end{aligned}
\]
Therefore,
\begin{equation}
  |R_\theta|\leq
    \left\{7\left(1+\frac{C_V}{\lambda_G}\right)
      +\frac\mu{\lambda_G}
        \left\|\frac{\nabla\bar\theta}{\bar\theta}\right\|_{L^\infty}^2\right\}
    \frac{\|q\|_{L^\infty}}\varepsilon Q_G(q).
  \label{eq:G-Rtheta-estimate}
\end{equation}

We now estimate the additional reaction remainder $R_G$.
Since $\|q\|_{L^\infty}\leq1/2$, for every $\tau\in[0,1]$ we have
\[
  \bar\theta+\tau z=\bar\theta(1+\tau q)
  \in\left[\frac12\bar\theta,\frac32\bar\theta\right].
\]
By Taylor's formula with integral remainder,
\[
\begin{aligned}
  &g(x,\bar\theta+z)-g(x,\bar\theta)-\partial_sg(x,\bar\theta)z
  =z^2\int_0^1(1-\tau)
    \partial_s^2g(x,\bar\theta+\tau z)\,d\tau.
\end{aligned}
\]
Substituting this into $R_G$ and using the definition of $M$ and $z=\bar\theta q$, we obtain
\[
\begin{aligned}
  |R_G|
  &\leq\frac1\varepsilon\int_{\Td}|z|^3
    \left\{\int_0^1(1-\tau)
      |\partial_s^2g(x,\bar\theta+\tau z)|\,d\tau\right\}\dd\\
  &\leq\frac{M}{2\varepsilon}\int_{\Td}\bar\theta^3|q|^3\dd\\
  &=\frac{M}{2\varepsilon}\int_{\Td}
    \bar\theta^2|q|\cdot\bar\theta q^2\dd\\
  &\leq\frac{M\|\bar\theta\|_{L^\infty}^2}{2\varepsilon}
    \|q\|_{L^\infty}\int_{\Td}\bar\theta q^2\dd\\
  &\leq\frac{M\|\bar\theta\|_{L^\infty}^2}{2\lambda_G}
    \frac{\|q\|_{L^\infty}}\varepsilon Q_G(q).
\end{aligned}
\]
The last inequality uses (H5$_G$).
Combining \eqref{eq:G-Rm-estimate}, \eqref{eq:G-Rtheta-estimate}, and the above estimate for $R_G$
with the definition of $C_R$, we obtain
\[
\begin{aligned}
  |R_m|+|R_\theta|+|R_G|
  &\leq C_R(\|r\|_{L^\infty}+\|q\|_{L^\infty})
       \{D_{\mathrm{flux}}+\mathcal{Q}_G\}\\
  &\leq\frac12\{D_{\mathrm{flux}}+\mathcal{Q}_G\}.
\end{aligned}
\]
The last inequality uses (H6$_G$).
Substituting this into \eqref{eq:G-energy-decomposition}
gives the dissipation estimate \eqref{eq:G-dissipation}.

Finally, we bound $|E|$ itself by the leading part.
Using $\int_{\Td}\bar m r\dd=0$ from \hypref{hyp:mass}
and the weighted Poincar\'e--Wirtinger inequality, we obtain
\[
  \int_{\Td}\bar m r^2\dd
  \leq C_P\int_{\Td}\bar m|\nabla r|^2\dd
  \leq C \mathcal{Q}_G(r,q).
\]
Thus $\nu\int_{\Td}n^2/\bar m\dd$ is bounded by $C \mathcal{Q}_G$.
Moreover, the inequality
$0\leq q-\log(1+q)\leq2q^2$ for $|q|\leq1/2$, together with (H5$_G$), gives
\[
  \frac1\varepsilon H_\theta
  \leq\frac2\varepsilon\int_{\Td}\bar\theta q^2\dd
  \leq\frac{2}{\varepsilon\lambda_G}Q_G(q)
  \leq C \mathcal{Q}_G(r,q).
\]
For the cross term, let $\baravg{v}$ denote the weighted mean of $v$ with respect to $\bar m$. Then
\[
\begin{aligned}
  \left|\int_{\Td}nv\dd\right|
  &=\left|\int_{\Td}\bar m r(v-\baravg{v})\dd\right|\\
  &\leq\sqrt{C_P}
    \left(\int_{\Td}\bar m r^2\dd\right)^{1/2}
    \left(\int_{\Td}\bar m|\nabla v|^2\dd\right)^{1/2}.
\end{aligned}
\]
Since $\nabla v=A+\nu\nabla r$,
\[
  \int_{\Td}\bar m|\nabla v|^2\dd
  \leq2D_{\mathrm{flux}}+2\nu^2\int_{\Td}\bar m|\nabla r|^2\dd
  \leq2\{D_{\mathrm{flux}}+\mathcal{Q}_G\}
\]
and therefore
\[
  \left|\int_{\Td}nv\dd\right|
  \leq C \mathcal{Q}_G^{1/2}
    (D_{\mathrm{flux}}+\mathcal{Q}_G)^{1/2}
  \leq C\{D_{\mathrm{flux}}+\mathcal{Q}_G\}.
\]
Combining the estimates for these three terms yields \eqref{eq:G-energy-control}.
Finally, \eqref{eq:G-dissipation} and \eqref{eq:G-energy-control} imply
\[
  -E'\geq\frac1{2C_{E,G}}|E|=\sigma_G|E|
\]
which gives \eqref{eq:G-signed-energy}.
\end{proof}

As in Proposition~\ref{prop:local-energy-di}, this result implies
\[
  |E(t)|\leq |E(0)|e^{-\sigma_Gt}
    +|E(T)|e^{-\sigma_G(T-t)}.
\]
The existence of solutions requires a separate analysis.

\par
\endgroup

\section{Sufficient conditions for (\ref{hyp:spectral-gap})}
Condition \eqref{hyp:spectral-gap} concerns the stationary solution. We now give a sufficient condition for it.
\begin{proposition}
   Let $\bar m$ and $\bar\theta$ be positive classical solutions of \eqref{eq:sKMR-FP} and \eqref{eq:sKMR-theta}, respectively.
   Set
    \begin{equation*}
        F(x):=K(x)-\varepsilon(\rho\ast\bar m)(x).
    \end{equation*}
    A sufficient condition for \eqref{hyp:spectral-gap} is
    \begin{equation*}
        \frac{3}{2}\inf_{x\in\Td}F(x)
        -
        \frac{1}{2}\sup_{x\in\Td}F(x)
        >0.
    \end{equation*}
\end{proposition}

\begin{proof}
    Let $\varphi\in H^1(\Td)$.
    Integration by parts on $\Td$, with no boundary contribution, gives
    \begin{align*}
        -\mu\int_{\Td}
        \varphi\nabla\bar\theta\cdot\nabla\varphi\,dx
        &=
        -\frac{\mu}{2}
        \int_{\Td}
        \nabla\bar\theta\cdot\nabla(\varphi^2)\,dx  \\
        &=
        \frac{\mu}{2}
        \int_{\Td}
        \Delta\bar\theta\,\varphi^2\,dx .
    \end{align*}
    Hence
    \begin{align*}
        Q_\theta(\varphi)
        &=
        \mu\int_{\Td}
        \bar\theta |\nabla\varphi|^2\,dx
        -
        \mu\int_{\Td}
        \varphi\nabla\bar\theta\cdot\nabla\varphi\,dx
        +
        \int_{\Td}
        \bar\theta^2\varphi^2\,dx  \\
        &=
        \mu\int_{\Td}
        \bar\theta |\nabla\varphi|^2\,dx
        +
        \int_{\Td}
        \left(
            \bar\theta^2+\frac{\mu}{2}\Delta\bar\theta
        \right)
        \varphi^2\,dx .
    \end{align*}
    Set
$
        V_\theta
        :=
        \bar\theta^2+\frac{\mu}{2}\Delta\bar\theta,
$.
Then \eqref{eq:sKMR-theta} gives
    \begin{align*}
        V_\theta
        &=
        \bar\theta^2
        +
        \frac{1}{2}
        \left(
            -\bar\theta F+\bar\theta^2
        \right)
        =
        \bar\theta
        \left(
            \frac{3}{2}\bar\theta
            -
            \frac{1}{2}F
        \right).
    \end{align*}

The classical maximum principle gives a lower bound for $\bar\theta$.
    Let $x_0\in\Td$ be a point where $\bar\theta$ attains its minimum.
    Then
$
        -\mu\Delta\bar\theta(x_0)\le 0
$.
Also, \eqref{eq:sKMR-theta} gives
$
        -\mu\Delta\bar\theta(x_0)
        =
        \bar\theta(x_0)
        \left(
            F(x_0)-\bar\theta(x_0)
        \right)
$.
Since $\bar\theta(x_0)>0$, it follows that
$
        F(x_0)-\bar\theta(x_0)\le 0
$.
Therefore
    \begin{equation*}
        \inf_{\Td}\bar\theta
        =
        \bar\theta(x_0)
        \ge
        F(x_0)
        \ge
        \inf_{\Td}F .
    \end{equation*}

    For every $x\in\Td$, we obtain
    \begin{align*}
        \frac{3}{2}\bar\theta(x)
        -
        \frac{1}{2}F(x)
        &\ge
        \frac{3}{2}\inf_{\Td}\bar\theta
        -
        \frac{1}{2}\sup_{\Td}F \\
        &\ge
        \frac{3}{2}\inf_{\Td}F
        -
        \frac{1}{2}\sup_{\Td}F \\
        &=
        \frac{3}{2}\inf_{x\in\Td}F(x)
        -
        \frac{1}{2}\sup_{x\in\Td}F(x) >0.
    \end{align*}
    The last step uses the assumption of the proposition.
    Thus
    \begin{equation*}
        V_\theta(x)
        =
        \bar\theta(x)
        \left(
            \frac{3}{2}\bar\theta(x)
            -
            \frac{1}{2}F(x)
        \right)
        \ge
        \left(\frac{3}{2}\inf_{x\in\Td}F(x)
        -
        \frac{1}{2}\sup_{x\in\Td}F(x)\right) \bar\theta(x)
        \qquad
        x\in\Td.
    \end{equation*}
    Consequently,
    \begin{align*}
        Q_\theta(\varphi)
        &=
        \mu\int_{\Td}
        \bar\theta |\nabla\varphi|^2\,dx
        +
        \int_{\Td}
        V_\theta\varphi^2\,dx \\
        &\ge
        \int_{\Td}
        V_\theta\varphi^2\,dx \\
        &\ge
        \left(\frac{3}{2}\inf_{x\in\Td}F(x)
        -
        \frac{1}{2}\sup_{x\in\Td}F(x)\right)
        \int_{\Td}
        \bar\theta\varphi^2\,dx .
    \end{align*}
    This is precisely (\ref{hyp:spectral-gap}).
\end{proof}
In particular, if $K$ is a positive constant, then (\ref{hyp:spectral-gap}) holds whenever the harvesting intensity $\varepsilon>0$ is sufficiently small.

\begin{corollary}
If $K=K_0>0$, a sufficient condition for \eqref{hyp:spectral-gap} is
\begin{equation*}
    2K_{0}> \varepsilon \left( 3 \sup_{\mathbb T^d} \rho \ast \bar{m}
    -\inf_{\mathbb T^d} \rho \ast \bar{m} \right).
\end{equation*}
\end{corollary}

\section{\texorpdfstring{$L^2$}{L2} turnpike estimates}
In this section, we assume
\hypref{hyp:even-kernel}--\hypref{hyp:smallness}
and derive
Proposition~\ref{prop:L2-turnpike_main} from Proposition~\ref{prop:local-energy-di}.
We first use
the differential inequality for the modified energy $E$
to bound its absolute value by exponential contributions from both endpoints of the time interval.
Integrating \eqref{eq:main-dissipation} in Proposition~\ref{prop:local-energy-di}
gives a bound on the time integral of $D_{\mathrm{flux}}+\mathcal{Q}$.
We then combine this bound with the differential inequalities for the individual components.
This yields $L^2$ turnpike estimates for $m-\bar m,\ \theta-\bar\theta$ and $\nabla u-\nabla \bar u$.

\begin{lemma}[Turnpike estimate for the energy]
Let $E\in C^1([0,T];\mathbb{R})$ satisfy
$  E'(t)\le -\sigma |E(t)|$ for all $t\in[0,T]$ and some $\sigma>0$.

Then, for every $t\in[0,T]$,
\begin{align}
  |E(t)|
  \le e^{-\sigma t}|E(0)|
  +e^{-\sigma(T-t)}|E(T)|
  \label{E-est}
\end{align}
holds.
\end{lemma}
\begin{proof}
By assumption,
\[
E'(t)\le -\sigma |E(t)|\le 0,
\qquad t\in[0,T]
\]
and hence $E$ is nonincreasing on $[0,T]$.
Fix $t\in[0,T]$ and consider the two possible signs of $E(t)$.

First suppose that $E(t)\ge 0$.
The monotonicity of $E$ gives, for every $s\in[0,t]$,
$E(s)\ge E(t)\ge 0$.
Therefore, on $[0,t]$,
\[
E'(s)\le -\sigma E(s)
\]
holds. Thus,
\[
\frac{d}{ds}\left(e^{\sigma s}E(s)\right)
=e^{\sigma s}\left(E'(s)+\sigma E(s)\right)
\le 0 .
\]
Integrating from $0$ to $t$, we obtain
\[
E(t)\le e^{-\sigma t}E(0).
\]
Since $E(0)\ge E(t)\ge0$ in this case,
\[
|E(t)|
\le e^{-\sigma t}|E(0)|.
\]

Similarly, if $E(t)<0$, we obtain
\[
|E(t)|
\le e^{-\sigma(T-t)}|E(T)|.
\]

Combining the two cases, we obtain, for every $t\in[0,T]$,
\[
|E(t)|
\le
e^{-\sigma t}|E(0)|
+
e^{-\sigma(T-t)}|E(T)|.
\]
\end{proof}

\subsection{\texorpdfstring{$L^2$}{L2} turnpike estimate for \texorpdfstring{$m$}{m}}
We first estimate
the harvester density $m$.
We derive a differential inequality by bounding the time derivative of the squared weighted $L^2$ norm
in terms of the dissipation $D_{\mathrm{flux}}+\mathcal{Q}$ in \eqref{eq:main-dissipation}
of Proposition~\ref{prop:local-energy-di}.
\begin{proof}[Proof of the estimate for $m$ in Proposition~\ref{prop:L2-turnpike_main}]
Since \(n=\bar m r\),
\[
\frac{d}{dt}\int_{\mathbb T^d}\bar m r^2 dx
=2\int_{\mathbb T^d} r\,\partial_t n\,dx .
\]
Substituting \eqref{eq:p2}, obtained from $\eqref{eq:KMR-FP}-\eqref{eq:sKMR-FP}$, namely
\[
\partial_t n=\operatorname{div}\{\bar m(\nu\nabla r-\nabla v)\}
-\operatorname{div}(n\nabla v)
\]
and integrating by parts, we obtain
\[
\frac{d}{dt}\int_{\mathbb T^d}\bar m r^2 dx
=
-2\nu\int_{\mathbb T^d}\bar m|\nabla r|^2dx
+2\int_{\mathbb T^d}\bar m\nabla r\cdot\nabla v\,dx
+2\int_{\mathbb T^d}n\nabla r\cdot\nabla v\,dx .
\]
Since
$A\mathrel{{=}}\nabla v-\nu\nabla r$,
we have \(\nabla v=A+\nu\nabla r\) and hence
\[
\frac{d}{dt}\int_{\mathbb T^d}\bar m r^2 dx
=
2\int_{\mathbb T^d}\bar m\nabla r\cdot A\,dx
+2\int_{\mathbb T^d}n\nabla r\cdot\nabla v\,dx .
\]
Young's inequality gives a constant $C>0$, independent of $t,T$, such that
\[
2\left|\bar m\nabla r\cdot A\right|
\le C\bar m\left(|\nabla r|^2+|A|^2\right).
\]
Also, \(|r|\le\delta\) implies \(|n|=\bar m|r|\le\delta\bar m\), and hence
\[
2|n\nabla r\cdot\nabla v|
\le 2\delta\bar m|\nabla r||\nabla v|
\le C\delta\bar m\left(|\nabla r|^2+|A|^2\right).
\]
Therefore, enlarging $C$ if necessary,
\[
\frac{d}{dt}\int_{\mathbb T^d}\bar m r^2 dx
\le
C\int_{\mathbb T^d}\bar m\left(|A|^2+|\nabla r|^2\right)dx .
\]
The definition of \(\mathcal{Q}\), together with
$Q_\theta\geq0$ from \eqref{hyp:spectral-gap},
gives a constant $C>0$, independent of $t,T$, such that
\[
D_{\mathrm{flux}}(t)+\mathcal{Q}(r(t),q(t))
\ge
C\int_{\mathbb T^d}\bar m\left(|A|^2+|\nabla r|^2\right)dx .
\]
Using \eqref{eq:main-dissipation} in Proposition~\ref{prop:local-energy-di}, namely
\[
-E'(t)\ge \frac{1}{2}\{D_{\mathrm{flux}}(t)+\mathcal{Q}(r(t),q(t))\}
\]
we obtain a constant $C>0$, independent of $t,T$, such that
\[
\frac{d}{dt}\int_{\mathbb T^d}\bar m r^2 dx
\le C(-E'(t)).
\]

Thus, for every \(0\le t_1<t_2\le T\),
\begin{align}
\int_{\mathbb T^d}\bar m r^2(t_2)dx
\le
\int_{\mathbb T^d}\bar m r^2(t_1)dx
+C\{E(t_1)-E(t_2)\}.
\label{eq:L^2_m_1}
\end{align}

Next,
we use the bound on $\int_{t-1}^{t}\|\sqrt{\bar m}r(s)\|_{L^2(\mathbb T^d)}^2\,ds$ obtained by integrating \eqref{eq:main-dissipation} in Proposition~\ref{prop:local-energy-di}
to estimate $\|\sqrt{\bar m}r(t)\|_{L^2(\mathbb T^d)}^2$ at each time.
For $t\geq 1$,
we choose a time $\tau$ at which the norm is bounded by its average over $[t-1,t]$.
We then use \eqref{eq:L^2_m_1} to extend this bound to time $t$.
For $0\leq t\leq\min\{1,T\}$,
instead of using the time average, we estimate $\|\sqrt{\bar m}r(s)\|_{L^2(\mathbb T^d)}^2$
by applying \eqref{eq:L^2_m_1} directly with $t_1=0$.

Let \(t\in[1,T]\).
By \eqref{eq:est_Q_by_P} and \eqref{eq:main-dissipation} in Proposition~\ref{prop:local-energy-di}, there exists $C>0$, independent of $t,T$, such that
\[
\int_{t-1}^{t}\int_{\mathbb T^d}\bar m r^2(s)dxds
\le
C\int_{t-1}^{t}\mathcal{Q}(r(s),q(s))ds
\le
C\{E(t-1)-E(t)\}.
\]
The mean value theorem for integrals then gives \(\tau\in(t-1,t)\) such that
\[
\int_{\mathbb T^d}\bar m r^2(\tau)dx
\le
C\{E(t-1)-E(t)\}.
\]
Applying \eqref{eq:L^2_m_1} with \(t_1=\tau,\ t_2=t\), we obtain
\[
\int_{\mathbb T^d}\bar m r^2(t)dx
\le
C\{E(t-1)-E(t)+E(\tau)-E(t)\}.
\]
Therefore,
\[
\int_{\mathbb T^d}\bar m r^2(t)dx
\le
C\left(|E(t-1)|+|E(\tau)|+|E(t)|\right).
\]
We use the energy turnpike estimate \eqref{E-est}. 
Since \(\tau\in(t-1,t)\),
we have
$e^{-\sigma \tau}\le e^\sigma e^{-\sigma t},
\
e^{-\sigma(T-\tau)}\le e^{-\sigma(T-t)}$.
Since $\sigma$ is independent of $T$, enlarging $C$ gives
\[
\int_{\mathbb T^d}\bar m r^2(t)dx
\le
C\left(e^{-\sigma t}|E(0)|
+e^{-\sigma(T-t)}|E(T)|\right).
\]

Finally, let \(t\in[0,{\min\{1,T\}}]\).
Applying \eqref{eq:L^2_m_1} with \(t_1=0,\ t_2=t\), we obtain
\[
\int_{\mathbb T^d}\bar m r^2(t)dx
\le
\int_{\mathbb T^d}\bar m r^2(0)dx
+C\{E(0)-E(t)\}.
\]
Therefore,
\[
\int_{\mathbb T^d}\bar m r^2(t)dx
\le
C\left(
\int_{\mathbb T^d}\bar m r^2(0)dx
+|E(0)|+|E(t)|
\right).
\]
We again use the energy turnpike estimate \eqref{E-est}. For \(t\in[0,{\min\{1,T\}}]\), \(e^{-\sigma t}\) is bounded below by a positive constant. Enlarging $C$ if necessary, we obtain
\[
\int_{\mathbb T^d}\bar m r^2(t)dx
\le
C\left(
\left(|E(0)|+\|\sqrt{\bar m}r(0)\|_{L^2(\mathbb T^d)}^2\right)e^{-\sigma t}
+
|E(T)|e^{-\sigma(T-t)}
\right).
\]

Thus, for every \(t\in[0,T]\),
\[
\|\sqrt{\bar m}r(t)\|_{L^2(\mathbb T^d)}^2
\le
\lambda_1 e^{-\sigma t}
+\lambda_2 e^{-\sigma(T-t)}
\]
holds,
where
\[
\lambda_1
=
C\left(|E(0)|+\|\sqrt{\bar m}r(0)\|_{L^2(\mathbb T^d)}^2\right),
\qquad
\lambda_2=C|E(T)|.
\]
\end{proof}

\subsection{\texorpdfstring{$L^2$}{L2} turnpike estimate for \texorpdfstring{$\theta$}{theta}}
We next estimate
the resource density $\theta$.
We pass from time averages to estimates at each time as in the argument for the harvester density $m$.

\begin{proof}[Proof of the estimate for $\theta$ in Proposition~\ref{prop:L2-turnpike_main}]
We first derive a differential inequality for
\[
\int_{\mathbb T^d}\bar\theta q^2\,dx
=
\int_{\mathbb T^d}
\frac{|\theta-\bar\theta|^2}{\bar\theta}\,dx.
\]

Writing \eqref{eq:p3}, obtained from $\eqref{eq:KMR-theta}\mathrel{{-}}\eqref{eq:sKMR-theta}$, in terms of $q$, we obtain
\[
\partial_t\theta
=
\mu\bar\theta\Delta {q}
+
2\mu\nabla\bar\theta\cdot\nabla {q}
-
\bar\theta^2{q(1+q)}
-
\varepsilon\bar\theta {(1+q)}(\rho*n).
\]

This gives
\[
\partial_t q
=
\mu\Delta q
+
2\mu\frac{\nabla\bar\theta}{\bar\theta}\cdot\nabla q
-
\bar\theta q(1+q)
-
\varepsilon(1+q)(\rho*n).
\]
Also, since
$\frac12{\frac{d}{dt}\int_{\mathbb T^d}\bar\theta q^2\,dx}
=
\int_{\mathbb T^d}\bar\theta q\,\partial_tq\,dx$,

integration by parts yields
\[
\begin{aligned}
\frac12{\frac{d}{dt}\int_{\mathbb T^d}\bar\theta q^2\,dx}
&=
-\mu\int_{\mathbb T^d}\bar\theta|\nabla q|^2\,dx
+
\mu\int_{\mathbb T^d}q\nabla\bar\theta\cdot\nabla q\,dx  \\
&\quad
-
\int_{\mathbb T^d}\bar\theta^2q^2(1+q)\,dx
-
\varepsilon\int_{\mathbb T^d}\bar\theta q(1+q)(\rho*n)\,dx .
\end{aligned}
\]
By \eqref{hyp:smallness}, \(\|q\|_{L^\infty}\le 1/2\), so \(1+q\ge 1/2\).
The third term is therefore nonpositive.
Thus,
\[
\frac12{\frac{d}{dt}\int_{\mathbb T^d}\bar\theta q^2\,dx}
\le
-\mu\int\bar\theta|\nabla q|^2
+
\mu\int q\nabla\bar\theta\cdot\nabla q
-
\varepsilon\int \bar\theta q(1+q)(\rho*n).
\]
Young's inequality gives
$\left|
q\nabla\bar\theta\cdot\nabla q
\right|
\le
\frac{1}{2}\bar\theta|\nabla q|^2
+\frac{1}{2}\frac{|\nabla \bar\theta|^2}{\bar\theta} q^2$.

Also,
\eqref{hyp:regularity} gives
$\bar\theta>0$ and boundedness of $|\nabla\bar\theta|$. Hence there exists $C>0$, independent of $t,T$, such that
\[
\left|
\mu\int q\nabla\bar\theta\cdot\nabla q
\right|
\le
\frac{\mu}{2}\int\bar\theta|\nabla q|^2
+
C\int\bar\theta q^2 .
\]
Next, by the boundedness of \(\rho,\bar\theta, 1+q\),
a similar calculation gives a constant $C>0$, independent of $t,T$, such that
\[
\left|
\varepsilon\int \bar\theta q(1+q)(\rho*n)
\right|
\le
C
\left(
\int\bar\theta q^2
+
\int \frac{n^2}{\bar m}
\right).
\]
By \eqref{hyp:spectral-gap} and the definition of \(\mathcal{Q}\), there exists $C>0$, independent of $t,T$, such that
\[
\int\bar\theta q^2\le C Q_\theta(q)\le C \mathcal{Q}(r,q),
\]
The weighted Poincar\'e--Wirtinger inequality and
$Q_\theta\geq0$
also give, after enlarging $C$ if necessary,
\[
\int\frac{n^2}{\bar m}
=
\int\bar m r^2
\le
C\int\bar m|\nabla r|^2
\le
C \mathcal{Q}(r,q).
\]
Thus, there exists $C>0$, independent of $t,T$, such that
\begin{align}
{\frac{d}{dt}\int_{\mathbb T^d}\bar\theta q^2\,dx}
\le
C \mathcal{Q}(r(t),q(t))
\le
C\{D_{\mathrm{flux}}(t)+\mathcal{Q}(r(t),q(t))\}.
\label{eq:theta_L^2_1}
\end{align}

Next, \eqref{eq:main-dissipation} in Proposition~\ref{prop:local-energy-di} gives, for every \({0\le}t_1<t_2{\le T}\),
\begin{align}
\frac{1}{2}\int_{t_1}^{t_2}
\{D_{\mathrm{flux}}(s)+\mathcal{Q}(r(s),q(s))\}\,ds
\le
E(t_1)-E(t_2).
\label{eq:theta_L^2_2}
\end{align}

Let \(t\in[1,T]\). By the mean value theorem for integrals,
there exists \(\tau\in(t-1,t)\) such that
\[
{\int_{\mathbb T^d}\bar\theta q^2(\tau)\,dx}
\le
\int_{t-1}^{t}{\int_{\mathbb T^d}\bar\theta q^2(s)\,dx}\,ds .
\]
Since \eqref{hyp:spectral-gap} gives \({\int_{\mathbb T^d}\bar\theta q^2(s)\,dx}\le C Q_\theta(q(s))\le C \mathcal{Q}(r(s),q(s))\),
\eqref{eq:theta_L^2_2} yields a constant $C>0$, independent of $t,T$, such that
\[
{\int_{\mathbb T^d}\bar\theta q^2(\tau)\,dx}
\le
C\{E(t-1)-E(t)\}.
\]
Integrating \eqref{eq:theta_L^2_1} over \([\tau,t]\) and using \eqref{eq:theta_L^2_2} again,
we obtain the following estimate after enlarging $C$.
Here we also use $E(\tau)\le E(t-1)$, which follows from the monotonicity of $E$.
\[
\begin{aligned}
{\int_{\mathbb T^d}\bar\theta q^2(t)\,dx}
&\le
{\int_{\mathbb T^d}\bar\theta q^2(\tau)\,dx}
+
C\int_{\tau}^{t}
\{D_{\mathrm{flux}}(s)+\mathcal{Q}(r(s),q(s))\}\,ds  \\
&\le
C\{E(t-1)-E(t)\}.
\end{aligned}
\]
Applying the energy turnpike estimate \eqref{E-est} and enlarging $C$, since $\sigma$ is independent of $T$, we obtain
\[
{\int_{\mathbb T^d}\bar\theta q^2(t)\,dx}
\le
C\left(
|E(0)|e^{-\sigma t}
+
|E(T)|e^{-\sigma(T-t)}
\right),
\qquad t\in[1,T].
\]

For \(t\in[0,{\min\{1,T\}}]\), we integrate \eqref{eq:theta_L^2_1} over \([0,t]\) to obtain
\[
{\int_{\mathbb T^d}\bar\theta q^2(t)\,dx}
\le
{\int_{\mathbb T^d}\bar\theta q^2(0)\,dx}
+
C\int_0^t
\{D_{\mathrm{flux}}(s)+\mathcal{Q}(r(s),q(s))\}\,ds .
\]
We use \eqref{eq:theta_L^2_2} and \eqref{E-est}, together with $1\le e^\sigma e^{-\sigma t}$ for $t\le1$. Enlarging $C>0$, which is independent of $t,T$, gives
\[
\begin{aligned}
{\int_{\mathbb T^d}\bar\theta q^2(t)\,dx}
\le & {\int_{\mathbb T^d}\bar\theta q^2(0)\,dx}+C\left(
|E(0)|
+
|E(t)|
\right)\\
\le &{\int_{\mathbb T^d}\bar\theta q^2(0)\,dx}+
C\left(
|E(0)|e^{-\sigma t}
+
|E(T)|e^{-\sigma(T-t)}
\right)\\
\le &
C\left(
\left\{|E(0)|+{\int_{\mathbb T^d}\bar\theta q^2(0)\,dx}\right\}e^{-\sigma t}
+
|E(T)|e^{-\sigma(T-t)}
\right).
\end{aligned}
\]
Combining these estimates, we obtain, for every \(t\in[0,T]\),
\[
\left\|
\frac{\theta(t)-\bar\theta}{\sqrt{\bar\theta}}
\right\|_{L^2}^2
\le
\lambda_3e^{-\sigma t}
+
\lambda_4e^{-\sigma(T-t)}.
\]

Here,
\[
 \lambda_3
    =C\left(
        |E(0)|+\left\|\frac{\theta_0-\bar\theta}{\sqrt{\bar\theta}}
        \right\|_{L^2(\mathbb{T}^d)}^2
    \right),
    \qquad\lambda_4
    =C|E(T)|.
\]

\end{proof}

\subsection{\texorpdfstring{$L^2$}{L2} turnpike estimate for the gradient of \texorpdfstring{$u$}{u}}
Finally,
we estimate the gradient of the value function $u$.
We use
the time-dependent density $m$ as a weight
to cancel the transport terms arising from the Fokker--Planck and Hamilton--Jacobi equations.
Set
\[
Y(t):=\frac12\int_{\mathbb T^d}m(t,x)|\nabla v(t,x)|^2\,dx.
\]

By \eqref{hyp:smallness}, $m$ and $\bar m$ are uniformly comparable.
Thus, a bound on this quantity yields
the desired estimate with weight $\bar m$.
\begin{proof}[Proof of the gradient estimate for $u$ in Proposition~\ref{prop:L2-turnpike_main}]
From $\eqref{eq:KMR-HJ}\mathrel{{-}}\eqref{eq:sKMR-HJ}$, we obtain
the equation for the Hamilton--Jacobi perturbation \eqref{eq:p1}, namely
\[
-\partial_t v
-\nu\Delta v
-{\nabla\bar u}\cdot {\nabla v}
={(\rho*z)}+\frac12|{\nabla v}|^2.
\]

Also,
substituting $\nabla u=\nabla\bar u+\nabla v$ into \eqref{eq:KMR-FP} gives
\[
\partial_t m
-\nu\Delta m
+\Div \left(m({\nabla\bar u}+{\nabla v})\right)
=0.
\]

Differentiating $Y$ with respect to time, we obtain
\[
Y'(t)
=
\frac12
\int_{\Td}
(\partial_t m)|{\nabla v}|^2\,dx
+
\int_{\Td}
m\,{\nabla v}\cdot\nabla(\partial_t v)\,dx.
\]

First, the equation \eqref{eq:KMR-FP} for $m$ and the periodic boundary conditions give
\[
\begin{aligned}
\frac12
\int_{\Td}
(\partial_t m)|{\nabla v}|^2\,dx
=\frac{\nu}{2}
\int_{\Td}
m\Delta(|{\nabla v}|^2)\,dx
+\frac12
\int_{\Td}
m({\nabla\bar u}+{\nabla v})\cdot\nabla(|{\nabla v}|^2)\,dx.
\end{aligned}
\]
Using
\[
\frac12\Delta(|{\nabla v}|^2)
=
|D^2v|^2
+
{\nabla v}\cdot\nabla\Delta v,
\qquad
\frac12\nabla(|{\nabla v}|^2)
=
D^2v\,{\nabla v},
\]
we obtain
\[
\begin{aligned}
\frac12
\int_{\Td}
(\partial_t m)|{\nabla v}|^2\,dx
={}&
\nu
\int_{\Td}
m|D^2v|^2\,dx
+\nu
\int_{\Td}
m\,{\nabla v}\cdot\nabla\Delta v\,dx\\
&+
\int_{\Td}
m\,{\nabla v}\cdot D^2v\,{\nabla\bar u}\,dx
+\int_{\Td}
m\,{\nabla v}\cdot D^2v\,{\nabla v}\,dx.
\end{aligned}
\]

Also,
differentiating \eqref{eq:p1} in space gives
\[
\nabla(\partial_t v)
=
-\nu\nabla\Delta v
-D^2v\,{\nabla\bar u}
-{D^2\bar u}{\nabla v}
-\nabla {(\rho*z)}
-D^2v\,{\nabla v}.
\]
Therefore,
\[
\begin{aligned}
\int_{\Td}
m\,{\nabla v}\cdot\nabla(\partial_t v)\,dx
={}&
-\nu
\int_{\Td}
m\,{\nabla v}\cdot\nabla\Delta v\,dx
-\int_{\Td}
m\,{\nabla v}\cdot D^2v\,{\nabla\bar u}\,dx
-\int_{\Td}
m\,{\nabla v}\cdot{D^2\bar u}{\nabla v}\,dx\\
&-
\int_{\Td}
m\,{\nabla v}\cdot\nabla {(\rho*z)}\,dx
-\int_{\Td}
m\,{\nabla v}\cdot D^2v\,{\nabla v}\,dx.
\end{aligned}
\]

Combining these identities, we obtain
\begin{align}
\label{u-L^2_1}
-Y'(t)
+\nu \int_{\Td} m|D^2v|^2\,dx
=\int_{\Td} m\,{\nabla v}\cdot{D^2\bar u}{\nabla v}\,dx
+\int_{\Td} m\,{\nabla v}\cdot\nabla {(\rho*z)}\,dx.
\end{align}

Since $\bar u\in C^2(\Td)$, there exists $C>0$, independent of $t,T$, such that
\[
\left|
\int_{\Td}
m\,{\nabla v}\cdot{D^2\bar u}{\nabla v}\,dx
\right|
\le
C Y(t).
\]

Next, since $\rho\in W^{1,1}(\Td)$,
Young's convolution inequality gives
\[
\|\nabla {(\rho*z)}\|_{L^2}
=
\|(\nabla\rho)*z\|_{L^2}
\le
\|\nabla\rho\|_{L^1}
\|z\|_{L^2}.
\]
Also, since
$z=\bar\theta q$ and
$0<{C_\theta^{-1}}\le\bar\theta\le C_\theta$,
there exists a constant $C>0$, independent of $t,T$, such that
\[
\|z\|_{L^2}^2
\le
C
\int_{\Td}
\bar\theta q^2\,dx.
\]
By \eqref{hyp:spectral-gap}, there exists $C>0$, independent of $t,T$, such that
$\int_{\Td}
\bar\theta q^2\,dx
\le
C Q_\theta(q)$.
Enlarging $C$ if necessary, we obtain
\[
\|\nabla {(\rho*z)}\|_{L^2}^2
\le
C Q_\theta(q)
\le
C \mathcal{Q}(r,q).
\]

By $m=\bar m(1+r)$ and \eqref{hyp:smallness},
\[
\frac12\bar m
\le
m
\le
\frac32\bar m.
\]
This and H\"{o}lder's inequality give a constant $C>0$, independent of $t,T$, such that
\[
\begin{aligned}
\left|
\int_{\Td}
m\,{\nabla v}\cdot\nabla {(\rho*z)}\,dx
\right|
&\le
\left(
\int_{\Td}
m|{\nabla v}|^2\,dx
\right)^{1/2}
\left(
\int_{\Td}
m|\nabla {(\rho*z)}|^2\,dx
\right)^{1/2}\\
&\le
C Y(t)^{1/2}
\mathcal{Q}(r,q)^{1/2}\\
&\le
C\left\{
Y(t)
+
\mathcal{Q}(r,q)
\right\}.
\end{aligned}
\]

The second term on the left-hand side of \eqref{u-L^2_1} is nonnegative. Thus, enlarging $C$, we obtain
\[
-Y'(t)
\le
C\left\{
Y(t)
+
\mathcal{Q}(r(t),q(t))
\right\}.
\]

Also, since
$\nabla v=A+\nu\nabla r$,
there exists a constant $C>0$, independent of $t,T$, such that
\[
\begin{aligned}
Y(t)
&\le
C
\int_{\Td}
\bar m|\nabla v|^2\,dx\\
&\le
C
\int_{\Td}
\bar m
\left(
|A|^2+\nu^2|\nabla r|^2
\right)\,dx\\
&\le
C
\left\{
D_{\mathrm{flux}}(t)
+
\mathcal{Q}(r(t),q(t))
\right\}.
\end{aligned}
\]
Hence, enlarging $C$ if necessary,
\[
-Y'(t)
\le
C
\left\{
D_{\mathrm{flux}}(t)
+
\mathcal{Q}(r(t),q(t))
\right\}.
\]

Combining this with the dissipation estimate \eqref{eq:main-dissipation} in Proposition~\ref{prop:local-energy-di}, we obtain a constant $C>0$, independent of $t,T$, such that
\[
-Y'(t)\le -C E'(t).
\]
Therefore, for every
$0\le t_1<t_2\le T$,

\begin{align}
\label{u-L^2_2}
Y(t_1)
\le
Y(t_2)
+
C\left\{
E(t_1)-E(t_2)
\right\}.
\end{align}

Next, let $t\in[0,T-1]$.
By the mean value theorem for integrals, there exists
$\tau\in(t,t+1)$ such that
\[
Y(\tau)
=
\int_t^{t+1}Y(s)\,ds.
\]
Also,
$Y(s)
\le
C
\left\{
D_{\mathrm{flux}}(s)
+
\mathcal{Q}(r(s),q(s))
\right\}$
and \eqref{eq:theta_L^2_2} give, after enlarging $C$ if necessary,
\[
Y(\tau)
\le
C\left\{
E(t)-E(t+1)
\right\}.
\]
Applying \eqref{u-L^2_2} with $t_1=t,\ t_2=\tau$, we obtain
\[
Y(t)
\le
C\left\{
|E(t)|
+
|E(t+1)|
+
|E(\tau)|
\right\}.
\]
For $s\in[t,t+1]$, we have $e^{-\sigma s}\le e^{-\sigma t}$ and $e^{-\sigma(T-s)}\le e^\sigma e^{-\sigma(T-t)}$.
Since $\sigma$ is independent of $T$, \eqref{E-est} gives the following bound after enlarging $C$.
\[
Y(t)
\le
C\left\{
e^{-\sigma t}|E(0)|
+
e^{-\sigma(T-t)}|E(T)|
\right\},
\qquad
t\in[0,T-1].
\]

Finally, let $t\in[{\max\{0,T-1\}},T]$.
Applying \eqref{u-L^2_2} with
$t_1=t$,
$t_2=T$,
we obtain
\[
Y(t)
\le
Y(T)
+
C\left\{
|E(t)|+|E(T)|
\right\}.
\]
For $t\in[{\max\{0,T-1\}},T]$, we have
$e^{-\sigma(T-t)}\ge e^{-\sigma}$.
Enlarging $C$ if necessary, we obtain
\[
Y(t)
\le
C e^{-\sigma t}|E(0)|
+
C e^{-\sigma(T-t)}
\left\{
|E(T)|+Y(T)
\right\}.
\]

Thus, enlarging $C$ once more, we obtain, for every $t\in[0,T]$,
\[
Y(t)
\le
C|E(0)|e^{-\sigma t}
+
C
\left\{
|E(T)|+Y(T)
\right\}
e^{-\sigma(T-t)}.
\]

Finally, since
\[
\frac{4}{3}Y(t)
\leq
\|\sqrt{\bar m}\nabla v(t)\|^2_{L^2(\Td)}
\leq 4Y(t),
\]

enlarging $C>0$, which is independent of $t,T$, gives
\[
\|\sqrt{\bar m}\nabla v(t)\|_{L^2(\Td)}^2
\le
C|E(0)|e^{-\sigma t}
+
C
\left\{
|E(T)|
+
\|\sqrt{\bar m}\nabla v(T)\|_{L^2(\Td)}^2
\right\}
e^{-\sigma(T-t)}.
\]

\end{proof}

This proves
all three estimates in Proposition~\ref{prop:L2-turnpike_main}.
The coefficients in these estimates
contain the endpoint values $E(0)$ and $E(T)$ of the modified energy.
In the next section,
we combine bounds on these endpoint energies with parabolic estimates to improve the constant in the $L^\infty$ turnpike estimate
and prove the existence of solutions near the stationary solution.

\section{Existence of solutions near a stationary solution}
In this section, we prove the existence of solutions near a stationary solution
(Proposition~\ref{prop_exsist}).
For each fixed $T>0$, we construct a classical solution satisfying the assumptions
when the initial and terminal data are sufficiently close to the stationary solution.

The proof uses the Leray--Schauder fixed point theorem.
In Section~5.1, we consider solutions satisfying an $L^\infty$ turnpike estimate
and use the $L^2$ turnpike estimate from Section~4
(Proposition~\ref{prop:L2-turnpike_main}) to improve its constant.
In Section~5.2, we construct a homotopy map for the system for the perturbations
around the stationary solution. The improved constant from Section~5.1 shows that the map
has no fixed points on the boundary of the prescribed neighborhood.
Together with the continuity and compactness of the map, this yields a fixed point.

{
We denote the deviations of the initial and terminal data from the stationary solution by
\[
v_T:=u_T-\bar u,\qquad
n_0:=m_0-\bar m,\qquad
z_0:=\theta_0-\bar\theta.
\]
The initial and terminal conditions for the original system for the perturbations \eqref{eq:p1}--\eqref{eq:p3} are
\[
v(T)=v_T,\qquad n(0)=n_0,\qquad z(0)=z_0.
\]
The initial datum for the relative perturbation $q=z/\bar\theta$ is determined by $z_0$ through
$
q_0:=\frac{z_0}{\bar\theta}
=\frac{\theta_0-\bar\theta}{\bar\theta}
$
and satisfies $q(0)=q_0$.
}

Using the constant $\sigma>0$ from Proposition~\ref{prop:local-energy-di}, set
\[
    \alpha:=\frac{\sigma}{d+1},
    \qquad
    \omega_\alpha(t)
    :=e^{-\alpha t}+e^{-\alpha(T-t)}.
\]
We aim to prove the following statement.
Let $(v,n,z)$ be a solution of the system for the perturbations \eqref{eq:p1}--\eqref{eq:p3} satisfying
\[
\|n(t)\|_{L^\infty}+\|q(t)\|_{L^\infty}
\le
\bar\delta
\omega_\alpha (t),
\qquad t\in[0,T].
\]
If the size of the initial and terminal data,
{
\[
\|n_0\|_{L^\infty}
+\|v_T\|_{W^{1,\infty}}
+\|z_0\|_{L^\infty},
\]
}
is sufficiently small, then
\[
\|n(t)\|_{L^\infty}+\|q(t)\|_{L^\infty}
\le
\frac{2}{3}\bar\delta
\omega_\alpha (t),
\qquad t\in[0,T].
\]
This improvement excludes fixed points on the boundary in Section~5.2.

\subsection{Improving the constant in the \texorpdfstring{\(L^\infty\)}{L-infty} turnpike estimate}
We first use the assumed $L^\infty$ turnpike estimate and the smallness of the terminal datum $v_T$
to obtain a bound on $\nabla v$ independent of $T$.
We then estimate the endpoint values of the modified energy in terms of the initial and terminal data
and substitute these bounds into the weighted $L^2$ estimates from Section~4.
Finally, parabolic estimates and interpolation inequalities allow us to reduce the constants
in the estimates for $n$ and $q$.

We begin with a gradient estimate for solutions of the perturbed Hamilton--Jacobi equation {\eqref{eq:p1}}.
We use the following estimate for linear backward parabolic equations.
\begin{lemma}[\cite{CN26}, Lemma A.4]
\label{lem:COCV-A4}
Let $z$ be the unique classical solution of the linear backward parabolic equation
\[
\begin{cases}
-\partial_t z-\nu\Delta z
=
\langle {\nabla}z,H\rangle+G,
& (t,x)\in(0,T)\times\Td,\\
z(T,\cdot)=z_T,
& x\in\Td.
\end{cases}
\]
Assume $H,G\in L^\infty([0,T]\times\Td)$ and
$z(t,\cdot)\in L^\infty(\Td),
\ t\in[0,T]$.
Then there exists a constant $C_H>0$, independent of $T$, such that
for every $t\in[0,T]$,
\[
\|{\nabla}z(t,\cdot)\|_{L^\infty(\Td)}
\le
2\|{\nabla}z_T\|_{L^\infty(\Td)}
+
C_H
\left(
\|G\|_{L^\infty([0,T]\times\Td)}
+
\|z\|_{L^\infty([0,T]\times\Td)}
\right).
\]
\end{lemma}

The perturbed Hamilton--Jacobi equation {\eqref{eq:p1}} contains a quadratic gradient term.
The Cole--Hopf transformation reduces it to a linear parabolic equation.
Combining this transformation with Lemma~\ref{lem:COCV-A4} gives the following uniform gradient estimate.
\begin{proposition}[Uniform estimate for the perturbed Hamilton--Jacobi equation]
\label{prop:HJ-smallness}
Assume {\eqref{hyp:even-kernel}}--{\eqref{hyp:smallness}}.
Let $v$ solve
\begin{equation}
\label{eq:v-perturbation}
    -\partial_t v
    -\nu\Delta v
    -\nabla\bar u\cdot\nabla v
    =
    \rho*(\bar\theta q)
    +\frac12|\nabla v|^2,
    \qquad
    v(T,\cdot)=v_T.
\end{equation}
There exist constants $\bar\delta>0$ and $\gamma_*>0$, independent of $T$, with the following property.
If $n,q$ and the terminal datum $v_T$ satisfy
\begin{equation}
\label{eq:bootstrap-nq}
    \|n(t)\|_{L^\infty(\mathbb T^d)}
    +
    \|q(t)\|_{L^\infty(\mathbb T^d)}
    \leq
    \bar\delta\,\omega_\alpha(t),
    \qquad t\in[0,T],
\end{equation}
and
\[
\|v_T\|_{W^{1,\infty}(\mathbb T^d)}
    \leq \gamma_*
\]
then
\begin{equation}
\label{eq:grad-v-small}
    \|\nabla v\|_{L^\infty([0,T]\times\mathbb T^d)}
    \leq 1.
\end{equation}
\end{proposition}

\begin{proof}
First, set
$C_\rho
    :=
    \|\rho\|_{L^1(\mathbb T^d)}
    \|\bar\theta\|_{L^\infty(\mathbb T^d)}$.
By \eqref{eq:bootstrap-nq},
\[
\begin{split}
    \|{\rho*(\bar\theta q(t))}\|_{L^\infty}
    &\leq
    \|\rho\|_{L^1}
    \|\bar\theta\|_{L^\infty}
    \|q(t)\|_{L^\infty}
    \\
    &\leq
    C_\rho\bar\delta\,\omega_\alpha(t),
\end{split}
\]

The comparison principle for \eqref{eq:v-perturbation}, applied with suitable subsolutions and supersolutions, gives
\[
    |v(t,x)|
    \leq
    \|v_T\|_{L^\infty}
    +
    \int_t^T\|{\rho*(\bar\theta q(s))}\|_{L^\infty}\,ds.
\]
Combining these estimates yields
\[
    \|v(t)\|_{L^\infty}
    \leq
    \|v_T\|_{L^\infty}
    +
    C_\rho\bar\delta
    \int_t^T\omega_\alpha(s)\,ds.
\]
Since
$\int_t^T\omega_\alpha(s)\,ds
    =\int_t^T e^{-\alpha s}\,ds
    +\int_t^T e^{-\alpha(T-s)}\,ds
    \leq\frac2\alpha$,
we obtain
\begin{equation}
\label{eq:Cv}
    \|v(t)\|_{L^\infty}
    \leq
    C_v,
    \qquad
    C_v
    :=\|v_T\|_{L^\infty}
    +\frac{2C_\rho}{\alpha}\bar\delta.
\end{equation}

Next, apply the Cole--Hopf transformation to \eqref{eq:v-perturbation} by setting
\begin{equation}
\label{eq:Cole-Hopf}
    w(t,x)
    :=
    \exp\left(
        \frac{ v(t,x)}{2\nu}
    \right).
\end{equation}
Then
\begin{equation}
\label{eq:w-linear}
    -\partial_t w
    -\nu\Delta w
    -\nabla\bar u\cdot\nabla w
    =\frac{{\rho*(\bar\theta q)}}{2\nu}w,
    \qquad
    w(T,x)
    =\exp\left(\frac{v_T(x)}{2\nu}\right).
\end{equation}
Setting $\hat w=w-1$, we find that $\hat w$ satisfies
\begin{equation}
\label{eq:hatw-linear}
    -\partial_t \hat w
    -\nu\Delta \hat w
    -\nabla\bar u\cdot\nabla \hat w
    =\frac{{\rho*(\bar\theta q)}}{2\nu}(\hat w+1),
    \qquad
    \hat w(T,x)
    =\exp\left(\frac{v_T(x)}{2\nu}\right)-1.
\end{equation}

Estimate \eqref{eq:Cv} gives
\begin{equation}
\label{eq:w-upper-lower}
    e^{-C_v/(2\nu)}-1
    \leq
    \hat w(t,x)
    \leq
    e^{C_v/(2\nu)}-1,
    \qquad
    \| \hat w\|_{L^\infty}
    \leq
    e^{C_v/(2\nu)}-1,
    \end{equation}
    \begin{equation}
      \|\nabla \hat w(T)\|_{L^\infty}
    \leq
    \frac{e^{\| v_T\|_{L^\infty}/(2\nu)}}{2\nu}
    \|\nabla v_T\|_{L^\infty},
    \qquad
     \left\|\frac{{(\rho*(\bar\theta q))}(1+\hat w)}{2\nu}\right\|_{L^\infty}
    \leq
    \frac{C_\rho}{\nu}\bar\delta e^{C_v/(2\nu)}.
\end{equation}

Applying Lemma~\ref{lem:COCV-A4} to \eqref{eq:hatw-linear},
we obtain a constant $C>0$, independent of $T$, such that
\begin{equation}
\label{eq:w-gradient}
    \sup_{t\in[0,T]}
    \|\nabla\hat w(t)\|_{L^\infty}
    \leq
    C\left(
        \|\nabla \hat w(T)\|_{L^\infty}
        +
        \left\|\frac{{(\rho*(\bar\theta q))}(1+\hat w)}{2\nu}\right\|_{L^\infty}
        +\|\hat w\|_{L^\infty}
    \right).
\end{equation}
Since
    $\nabla v
    =2\nu\frac{\nabla \hat w}{w}$,
we obtain
\begin{equation}
\label{eq:final-gradient-estimate}
    \|\nabla v\|_{L^\infty([0,T]\times\mathbb T^d)}
    \leq
    C e^{C_v/\nu}
    \left(
        e^{\| v_T\|_{L^\infty}/(2\nu)}
    \|\nabla v_T\|_{L^\infty}
        +\bar\delta
        +\left|e^{C_v/(2\nu)}-1\right|
    \right)
\end{equation}
with a constant $C$ independent of $T$.

Finally, \eqref{eq:Cv} implies that $C_v\to0$ as $\|v_T\|_{W^{1,\infty}}+\bar\delta\to0$.
Thus, the right-hand side of \eqref{eq:final-gradient-estimate} also tends to zero.
Choosing $\gamma_*>0$ and $\bar\delta>0$ sufficiently small makes the right-hand side of \eqref{eq:final-gradient-estimate} at most one
whenever $\|v_T\|_{W^{1,\infty}}\le\gamma_*$. This proves the proposition.
\end{proof}

The constants in the $L^2$ estimates of Proposition~\ref{prop:L2-turnpike_main}
involve the endpoint values $E(0)$ and $E(T)$ of the modified energy.
These also depend on unknown endpoint values of the solution.
To express the estimates in terms of small initial and terminal data,
we therefore need to control these endpoint energies.
The gradient estimate above and monotonicity of the modified energy yield the following bound.

\begin{lemma}[Endpoint estimate for the modified energy]
\label{lem:E-estimate}
Assume {\eqref{hyp:even-kernel}}--{\eqref{hyp:smallness}} and
$\|\nabla v(0)\|_{L^\infty(\mathbb T^d)}\le 1$.
Set
$X_0:=
\left\|
\frac{{n_0}}{\sqrt{\bar m}}
\right\|_{L^2(\mathbb T^d)}$.
Then there exists a constant $C>0$ such that
\[
|E(0)|+|E(T)|
\le
C\left\{
X_0+\nu X_0^2
+\frac1\varepsilon H_\theta(0)
+
\frac{C_P}{2\nu}
\left\|
\sqrt{\bar m}\,\nabla {v_T}
\right\|_{L^2(\mathbb T^d)}^2
\right\}.
\]
\end{lemma}

\begin{proof}
By the definition of the modified energy,
\[
E(0)
=
-\int_{\mathbb T^d}{n_0}v(0)\,dx
+
\nu\int_{\mathbb T^d}\frac{{n_0}^2}{\bar m}\,dx
+
\frac1\varepsilon H_\theta(0).
\]

We denote the weighted mean by
$\langle v(t)\rangle_{\bar m}
{=}\int_{\mathbb T^d}\bar m\,v(t)\,dx$.
Then H\"older's inequality gives
\begin{align}
\left|\int_{\mathbb T^d}{n_0}v(0)\,dx\right|
&=
\left|
\int_{\mathbb T^d}
{n_0}\left(v(0)-\langle v(0)\rangle_{\bar m}\right)\,dx
\right|\\
&=
\left|
\int_{\mathbb T^d}
\frac{{n_0}}{\sqrt{\bar m}}\,
\sqrt{\bar m}
\left(v(0)-\langle v(0)\rangle_{\bar m}\right)\,dx
\right|\\
&\le
\left\|\frac{{n_0}}{\sqrt{\bar m}}\right\|_{L^2}
\left\|
\sqrt{\bar m}
\left(v(0)-\langle v(0)\rangle_{\bar m}\right)
\right\|_{L^2}\\
&=
X_0
\left(
\int_{\mathbb T^d}
\bar m
\left|v(0)-\langle v(0)\rangle_{\bar m}\right|^2dx
\right)^{1/2}.
\end{align}
The assumed bound on
$\nabla v(0)$
and the weighted Poincar\'e--Wirtinger inequality
$\int_{\mathbb T^d}
\bar m|v-\langle v\rangle_{\bar m}|^2\,dx
\le
C_P\int_{\mathbb T^d}\bar m|\nabla v|^2\,dx
$
give, with a constant $C>0$ independent of $T$, enlarged if necessary,
\[
\left|\int_{\mathbb T^d}{n_0}v(0)\,dx\right|
\le
C X_0\|\sqrt{\bar m}\nabla v(0)\|_{L^2}
\le
{C} X_0
\]

\[
E(0)
\le
{C}X_0+\nu X_0^2+\frac1\varepsilon H_\theta(0).
\]

At the terminal time, the definition of the modified energy gives
\[
E(T)
=
-\int_{\mathbb T^d}n(T){v_T}\,dx
+
\nu\int_{\mathbb T^d}\frac{n(T)^2}{\bar m}\,dx
+
\frac1\varepsilon H_\theta(T).
\]
Since
\(
\int_{\mathbb T^d}n(T)\,dx=0
\),
we have
\[
\int_{\mathbb T^d}n(T){v_T}\,dx
=
\int_{\mathbb T^d}
n(T)
\left(
{v_T}-\langle {v_T}\rangle_{\bar m}
\right)
dx.
\]
Young's inequality therefore yields
\[
\begin{aligned}
-\int_{\mathbb T^d}n(T){v_T}\,dx
\ge&
-\frac{\nu}{2}
\int_{\mathbb T^d}
\frac{n(T)^2}{\bar m}\,dx
-\frac1{2\nu}
\int_{\mathbb T^d}
\bar m
\left|
{v_T}-\langle {v_T}\rangle_{\bar m}
\right|^2dx.
\end{aligned}
\]
Since
\(
H_\theta(T)\ge0
\),
we also have
\[
\begin{aligned}
E(T)
\ge&
\frac{\nu}{2}
\int_{\mathbb T^d}
\frac{n(T)^2}{\bar m}\,dx-\frac1{2\nu}
\int_{\mathbb T^d}
\bar m
\left|
{v_T}-\langle {v_T}\rangle_{\bar m}
\right|^2dx\\
\ge&
-\frac1{2\nu}
\int_{\mathbb T^d}
\bar m
\left|
{v_T}-\langle {v_T}\rangle_{\bar m}
\right|^2dx.
\end{aligned}
\]
The weighted Poincar\'e--Wirtinger inequality
$\int_{\Td} \bar m\,|v-\langle v\rangle_{\bar m}|^2
\le
C_P
\int_{\Td} \bar m\,|\nabla v|^2$
then implies
\[
E(T)
\ge
-
\frac{C_P}{2\nu}
\left\|
\sqrt{\bar m}\,\nabla {v_T}
\right\|_{L^2(\mathbb T^d)}^2.
\]
We have thus obtained
\[
\begin{aligned}
E(0)
&\le
{C}X_0+\nu X_0^2+\frac1\varepsilon H_\theta(0),\\
E(T)
&\ge
-
\frac{C_P}{2\nu}
\left\|
\sqrt{\bar m}\,\nabla {v_T}
\right\|_{L^2(\mathbb T^d)}^2.
\end{aligned}
\]
Since $E'(t)\le -\sigma |E(t)|$, we have
$E(T)\le E(0)$
and hence, enlarging $C$ if necessary,
\[
|E(0)|+|E(T)|
\le
C\left\{
X_0+\nu X_0^2
+\frac1\varepsilon H_\theta(0)
+
\frac{C_P}{2\nu}
\left\|
\sqrt{\bar m}\,\nabla {v_T}
\right\|_{L^2(\mathbb T^d)}^2
\right\}.
\]
\end{proof}

Substituting the endpoint estimate from Lemma~\ref{lem:E-estimate}
into Proposition~\ref{prop:L2-turnpike_main} controls the constants in the $L^2$ turnpike estimates
by the initial and terminal data. In particular, these constants tend to zero
as the data approach the stationary solution. We record these estimates in the following lemma.
\begin{lemma}[Weighted $L^2$ turnpike estimates]
\label{lem:L2-turnpike-}
Assume {\eqref{hyp:even-kernel}}--{\eqref{hyp:smallness}} and
\[
\|{\nabla u(0)-\nabla\bar u}\|_{L^\infty(\mathbb T^d)}\le 1.
\]
Then there exist constants independent of $T$, denoted by
$\Lambda_1,\Lambda_2,\Lambda_3,\Lambda_4,
\Lambda_5,\Lambda_6\geq0$
such that, for every $t\in[0,T]$,
\[
\left\|
\frac{m(t)-\bar m}{\sqrt{\bar m}}
\right\|_{L^2(\mathbb T^d)}^2
\le
\Lambda_1e^{-\sigma t}
+
\Lambda_2e^{-\sigma(T-t)},
\]
\[
\left\|
\frac{\theta(t)-\bar\theta}{\sqrt{\bar\theta}}
\right\|_{L^2(\mathbb T^d)}^2
\le
\Lambda_3e^{-\sigma t}
+
\Lambda_4e^{-\sigma(T-t)},
\]
and
\[
\left\|
\sqrt{\bar m}\,{(\nabla u(t)-\nabla\bar u)}
\right\|_{L^2(\mathbb T^d)}^2
\le
\Lambda_5e^{-\sigma t}
+
\Lambda_6e^{-\sigma(T-t)}
\]
hold. The constants $\Lambda_i$ can be chosen as
\[
\begin{aligned}
\Lambda_1
&=
C\left(
X_0+(\nu+1) X_0^2
+\frac1\varepsilon H_\theta(0)
+
\frac{C_P}{2\nu}
X_1^2
\right),\\
\Lambda_2
&=C\left(
X_0+\nu X_0^2
+\frac1\varepsilon H_\theta(0)
+
\frac{C_P}{2\nu}
X_1^2
\right),
\end{aligned}
\]
\[
\begin{aligned}
\Lambda_3
&=
C\left(
X_0+\nu X_0^2
+\frac1\varepsilon H_\theta(0)
+
\frac{C_P}{2\nu}
X_1^2
+
\left\|
\frac{\theta_0-\bar\theta}{\sqrt{\bar\theta}}
\right\|_{L^2(\mathbb T^d)}^2
\right),\\
\Lambda_4
&=C\left(
X_0+\nu X_0^2
+\frac1\varepsilon H_\theta(0)
+
\frac{C_P}{2\nu}
X_1^2
\right),
\end{aligned}
\]
\[
\begin{aligned}
\Lambda_5
&=C\left(
X_0+\nu X_0^2
+\frac1\varepsilon H_\theta(0)
+
\frac{C_P}{2\nu}
X_1^2
\right),\\
\Lambda_6
&=
C\left(
X_0+\nu X_0^2
+\frac1\varepsilon H_\theta(0)
+
(\frac{C_P}{2\nu}+1)
X_1^2
\right)
\end{aligned}
\]
where
\[
X_0:=
\left\|
\frac{{n_0}}{\sqrt{\bar m}}
\right\|_{L^2(\mathbb T^d)},
\qquad
X_1:=\left\|
\sqrt{\bar m}\,{(\nabla u_T-\nabla\bar u)}
\right\|_{L^2(\mathbb T^d)}.
\]
\end{lemma}

We next prepare a parabolic a priori estimate to pass from the $L^2$ bounds to $L^\infty$ bounds.
To apply it to both the $n$ and $q$ equations, we include an inhomogeneous term $H$
in addition to the forcing term in divergence form.
The proof follows \cite[Lemma A.2]{CN26}.
The additional contribution of $H$ is handled by heat semigroup estimates,
just as the zeroth-order term arising from the time cutoff in that proof.
We therefore omit the proof.
\begin{theorem}[A priori estimate]
\label{thm:linfty-estimate-general-diffusion}
Let $T>0$, $C>0$, and assume
\[
F\in L^\infty\left([0,T]\times\Td\right),
\qquad
G\in L^\infty\left([0,T];L^{d+1}(\Td)\right),
\qquad
H\in L^\infty\left([0,T];L^{d}(\Td)\right).
\]
Let $z$ satisfy
\begin{equation}
\label{eq:general-divergence-equation}
\begin{cases}
\partial_t z-C\Delta z
 =-\operatorname{div}(zF)
-\operatorname{div}G
+H,
&
(t,x)\in(0,T)\times\Td,
\\
z(0,x)=z_{0}(x),
&
x\in\Td.
\end{cases}
\end{equation}
Let $c_d>0$ depend only on the dimension $d$ and set
{
\[
\tau=
\min\left\{
1,\,
\frac{C}{16c_d^2(1+\|F\|_{L^\infty((0,T)\times\Td)}^2)}
\right\}.
\]
}
Then there exists a constant
\[
K=K\left(
d,C,\tau,
\|F\|_{L^\infty((0,T)\times\Td)}
\right)>0
\]
such that, for every $t\in{(0,T]}$,
\[
\begin{aligned}
\|z(t)\|_{L^\infty}
\le {}&
2(\mathbf 1_{[0,\tau]}(t)\|z_0\|_{L^\infty})
\\
&+K\left(
\|G\|_{L^\infty({(\max\{0,t-\tau\},t)};L^{d+1})}
+\|H\|_{L^\infty({(\max\{0,t-\tau\},t)};L^{d})}\right.
\\
&\qquad\left.+\mathbf 1_{[\tau,T]}(t)\|z\|_{L^\infty({(\max\{0,t-\tau\},t)};L^d)}
\right).
\end{aligned}
\]
\end{theorem}

We now prove the improvement of the constant in the turnpike estimate.
Combining the $L^2$ bounds of Lemma~\ref{lem:L2-turnpike-} with interpolation inequalities
and the a priori estimate in Theorem~\ref{thm:linfty-estimate-general-diffusion}
gives $L^\infty$ estimates for $n$ and $q$.
Choosing the initial and terminal data sufficiently small makes the constant
smaller than the one assumed initially.
\begin{proposition}[Improvement of the $L^\infty$ turnpike estimate]
\label{prop:n-Linfty-improvement}
Assume {\eqref{hyp:even-kernel}}--{\eqref{hyp:smallness}} and
\[
\|n(t)\|_{L^\infty}+\|q(t)\|_{L^\infty}
\le
\bar\delta
\omega_\alpha (t),
\qquad t\in[0,T],
\]
Then there exists $\bar\gamma\in(0,\gamma_*)$ such that, if
\[
\|{n_0}\|_{L^\infty}
+\|u_T-\bar u\|_{W^{1,\infty}}
+\left\|z_0\right\|_{L^\infty}
<\bar\gamma,
\]
then
\[
\|n(t)\|_{L^\infty}+\|q(t)\|_{L^\infty}
\le
\frac{2\bar\delta}{3}\omega_\alpha(t),
\qquad t\in[0,T],
\]
Here $\bar\delta, \gamma_*$ are as in Proposition~{\ref{prop:HJ-smallness}}.
\end{proposition}

\begin{proof}
We estimate the solution for $t\in(0,T]$.
The assumptions allow us to apply Proposition~{\ref{prop:HJ-smallness}}, which gives
    $\|\nabla v\|_{L^\infty([0,T]\times\mathbb T^d)}
    \leq 1$.
Using the zero-flux identity from Lemma~\ref{lem:fluz},
we recall that the perturbation equation {\eqref{eq:p2}} for $n$ can be written as
\[
\partial_tn-\nu\Delta n
+
\operatorname{div}
\left(
n(\nabla\bar u+\nabla v)
\right)
=
-\operatorname{div}(\bar m\nabla v).
\]
Apply Theorem~\ref{thm:linfty-estimate-general-diffusion} with
$z=n,\
F=\nabla\bar u+\nabla v,\
G=\bar m\nabla v$.
Since $\|\nabla v\|_{L^\infty}\le 1$,
the norm $\|F\|_{L^\infty}$ is bounded by a constant independent of $T$.
Thus there exist $\tau>0$ and $C>0$, independent of $T$, such that
\[
\begin{aligned}
\|n(t)\|_{L^\infty}
\le {}&
2\mathbf 1_{[0,\tau]}(t)\|{n_0}\|_{L^\infty}
\\
&+C\|\bar m\nabla v\|_{L^\infty({(\max\{0,t-\tau\},t)};L^{d+1})}
\\
&+C\|n\|_{L^\infty({(\max\{0,t-\tau\},t)};L^d)}.
\end{aligned}
\]
We first estimate the term involving $\nabla v$.
Boundedness of $\bar m$ and the interpolation inequality \cite[Chapter 4, Remark 2]{B11} give
\[
\begin{aligned}
\|\bar m\nabla v\|_
{L^\infty({(\max\{0,t-\tau\},t)};L^{d+1})}
\le&
C\|\nabla v(s)\|_{L^\infty({(\max\{0,t-\tau\},t)};L^{d+1})}\\
\le&
C\|\nabla v(s)\|_{L^\infty([0,T]\times \Td)}^{1-\frac{2}{d+1}}
\sup_{s\in {(\max\{0,t-\tau\},t)}}
\|\nabla v(s)\|_{L^2}^{\frac{2}{d+1}}.
\end{aligned}
\]
Lemma~\ref{lem:L2-turnpike-} implies
\[
\|\nabla v(s)\|_{L^2}^{\frac{2}{d+1}}
\le
C
\left(
\Lambda_5^{\frac1{d+1}}
e^{-\frac{\sigma}{d+1}s}
+
\Lambda_6^{\frac1{d+1}}
e^{-\frac{\sigma}{d+1}(T-s)}
\right).
\]
Next, we estimate the $L^d$ norm of $n$.
For $d>2$, the interpolation inequality \cite[Chapter 4, Remark 2]{B11} gives
\[
\|n(s)\|_{L^d}
\le
\|n(s)\|_{L^\infty}^{1-\frac2d}
\|n(s)\|_{L^2}^{\frac2d}.
\]
Using the assumptions and Lemma~\ref{lem:L2-turnpike-}, we obtain
\[
\begin{aligned}
\|n(s)\|_{L^d}
\le
C\bar\delta^{1-\frac2d}
\omega_\alpha(s)^{1-\frac2d}
\left(
\Lambda_1e^{-\sigma s}
+
\Lambda_2e^{-\sigma(T-s)}
\right)^{\frac1d}.
\end{aligned}
\]
Consequently,
\[
\|n(s)\|_{L^d}
\le
C\bar\delta^{1-\frac2d}
\left(
\Lambda_1^{\frac1d}e^{-\beta s}
+
\Lambda_2^{\frac1d}e^{-\beta(T-s)}
\right),
\]
where
\[
\beta
:=
\alpha\left(1-\frac2d\right)
+
\frac{\sigma}{d}.
\]
By the definition of $\alpha$,
\[
\beta-\alpha
=
\frac{\sigma-2\alpha}{d}
\ge0.
\]
Thus all exponential terms on the right-hand side are bounded by $\omega_\alpha$.

{
For $d\le 2$, the finite measure of $\mathbb T^d$ gives
\[
\|n(s)\|_{L^d(\mathbb T^d)}
\le C\|n(s)\|_{L^2(\mathbb T^d)}.
\]
Using the weighted $L^2$ turnpike estimate and the upper and lower bounds on $\bar m$, we obtain
\[
\|n(s)\|_{L^d}
\le
C\left(
\Lambda_1^{1/2}e^{-\sigma s/2}
+
\Lambda_2^{1/2}e^{-\sigma(T-s)/2}
\right).
\]
Since $\alpha=\sigma/(d+1)\le \sigma/2$, we have
\[
e^{-\sigma s/2}+e^{-\sigma(T-s)/2}
\le
\omega_\alpha(s)
\]
and hence
\[
\|n(s)\|_{L^d}
\le
C\left(\Lambda_1^{1/2}+\Lambda_2^{1/2}\right)
\omega_\alpha(s).
\]
}

Moreover, there exists a constant $C_{\alpha,\tau}>0$, depending only on $\alpha,\tau$, such that, if $s\in{(\max\{0,t-\tau\},t)}$, then
\[
\omega_\alpha(s)
\le
C_{\alpha,\tau}\omega_\alpha(t).
\]
Thus,
{
adjusting the constant $C>0$, we obtain
\[
\begin{aligned}
\|n(t)\|_{L^\infty}
\le
2\mathbf 1_{[0,\tau]}(t)\|{n_0}\|_{L^\infty}
+
C\left[
\Lambda_5^{\frac1{d+1}}
+
\Lambda_6^{\frac1{d+1}}
+
\Lambda_1^{\frac{1}{\max\{2,d\}}}
+
\Lambda_2^{\frac{1}{\max\{2,d\}}}
\right]
\omega_\alpha(t).
\end{aligned}
\]
}

By Lemma~\ref{lem:L2-turnpike-}, as
$\|{n_0}\|_{L^\infty}
+\|u_T-\bar u\|_{W^{1,\infty}}
+\left\|z_0\right\|_{L^\infty}
{\to}0$, we have $\Lambda_i{\to}0,
\ (i=1,\dots,6)$.
Moreover,
\[
1_{[0,\tau]}(t)\|{n_0}\|_\infty
\le
e^{\alpha\tau}\|{n_0}\|_\infty\omega_\alpha(t)
\]
and, as $\bar\gamma\to 0$, we have
$
\|{n_0}\|_{L^\infty}{\to} 0
$.
Thus, choosing $\bar\gamma>0$ sufficiently small, we can ensure that
{
\[
2\mathbf 1_{[0,\tau]}(t)\|{n_0}\|_{L^\infty}
+
C\left[
\Lambda_5^{\frac1{d+1}}
+
\Lambda_6^{\frac1{d+1}}
+
\Lambda_1^{\frac{1}{\max\{2,d\}}}
+
\Lambda_2^{\frac{1}{\max\{2,d\}}}
\right]
\omega_\alpha(t)
\le
\frac{\bar\delta}{3}\omega_\alpha(t),\qquad t\in[0,T].
\]
}
It follows that
\[
\|n(t)\|_{L^\infty}
\le
\frac{\bar\delta}{3}
\left(
e^{-\alpha t}
+
e^{-\alpha(T-t)}
\right).
\]

Next, recall that the equation for $q$ in {\eqref{eq:p3}} is
\[
\partial_tq-\mu\Delta q
=
\operatorname{div}
\left(
2\mu q\frac{\nabla\bar \theta}{\bar \theta}
\right)
-
\left[
\operatorname{div}
\left(
2\mu \frac{\nabla\bar \theta}{\bar \theta}
\right)
+\bar \theta(1+q)
\right]q
-\varepsilon(\rho\ast n)(1+q).
\]
Apply Theorem~\ref{thm:linfty-estimate-general-diffusion} with
\[
z=q,\qquad
F=-2\mu \frac{\nabla\bar \theta}{\bar \theta},\qquad
G=0,\qquad
H=-\left[
\operatorname{div}
\left(
2\mu \frac{\nabla\bar \theta}{\bar \theta}
\right)
+\bar \theta(1+q)
\right]q
-\varepsilon(\rho\ast n)(1+q).
\]
By assumption~\eqref{hyp:regularity},
the norm $\|F\|_{L^\infty}$ is bounded by a constant independent of $T$.
Thus there exist $\tau>0$ and $C>0$, independent of $T$, such that
\[
\begin{aligned}
\|q(t)\|_{L^\infty}
&\le
2\mathbf 1_{[0,\tau]}(t)\|q_0\|_{L^\infty}
\\
&\quad+C\|\operatorname{div}
\left(
2\mu \frac{\nabla\bar \theta}{\bar \theta}
\right)
+\bar \theta(1+q)
\|_{L^\infty({(\max\{0,t-\tau\},t)};L^{\infty})}
\|q\|_{L^\infty({(\max\{0,t-\tau\},t)};L^d)}
\\
&\quad+C\|\varepsilon(\rho\ast n)(1+q)\|_
{L^\infty({(\max\{0,t-\tau\},t)};L^d)}
+C\|q\|_{L^\infty({(\max\{0,t-\tau\},t)};L^d)}\\.
\end{aligned}
\]
Since $q, \bar\theta, \rho$ are bounded independently of $T$, we obtain
\[
\begin{aligned}
\|q(t)\|_{L^\infty}
\le {}&
2\mathbf 1_{[0,\tau]}(t)\|q_0\|_{L^\infty}
\\
&+C'\|q\|_{L^\infty({(\max\{0,t-\tau\},t)};L^d)}
+C'\|n\|_
{L^\infty({(\max\{0,t-\tau\},t)};L^d)}
\\
&+C\|q\|_{L^\infty({(\max\{0,t-\tau\},t)};L^d)}.
\end{aligned}
\]
Combining the interpolation inequality \cite[Chapter 4, Remark 2]{B11} with the $L^2$ turnpike estimate,
as in the argument for $n$, proves the conclusion.
\end{proof}

\subsection{Existence of solutions satisfying (\ref{hyp:smallness}) near a stationary solution}
We prove existence by the Leray--Schauder fixed point theorem.
We first construct a suitable homotopy map ${\mathcal F}$.

Let $(\bar u,\bar m,\bar\theta,\lambda)$ be a stationary solution and set
$
  U^*(t,x):=\bar u(x)+\lambda(T-t)
$.
Define the perturbations around this stationary solution by
\[
  v:=u-U^*,\qquad
  n:=m-\bar m,\qquad
  q:=\frac{\theta-\bar\theta}{\bar\theta}.
\]
Here $v_T,n_0,q_0$ are the initial and terminal data defined at the beginning of Section~5.
Equip the space
$
X=C([0,T]\times \Td) \times C([0,T]\times \Td)
$
with the norm
\[
\|(\eta, \xi)\|_{X}
=\sup_{t \in [0,T]} \frac{\|\eta{(t,\cdot)}\|_{L^\infty{(\Td)}}+
             \|\xi{(t,\cdot)}\|_{L^\infty{(\Td)}}}{e^{-\alpha t}+e^{-\alpha (T-t)}}.
\]
Then $X$ is a Banach space.

Let $B_\delta^X(0)$ denote the open ball in $X$ centered at the origin with radius $\delta$.
For a homotopy parameter $s\in[0,1]$, define the map
\[
  \mathcal F{:\overline{B_\delta^X(0)}\times[0,1]\mapsto X}
\]
by specifying $\mathcal F((\eta,\xi),s)=(\tilde n,\tilde q)$ in three steps.
First, given the input $\xi$, solve
\[
\begin{cases}
 -\partial_t
 v-\nu\Delta v-\nabla\bar u\cdot\nabla v
 =
 s\,\rho*(\bar\theta\xi)+\dfrac12|\nabla v|^2,
\\
 v(T,x)=s\,v_T(x)
\end{cases}
\]
as in
\cite{Pazy1983}.
Next, using the unique classical solution $v$, solve
\[
\begin{cases}
 \partial_t\tilde n
 =
 \operatorname{div}\left\{\bar m\left(\nu\nabla{\left(\frac{\tilde n}{\bar m}\right)}-\nabla v\right)\right\}
 -\operatorname{div}(\tilde n\nabla v),
\\
 \tilde n(0,x)=s\,n_0(x).
\end{cases}
\]
Setting $U=\bar m+\tilde n$ and $V=\bar u+v$, we can rewrite this initial-value problem as
\begin{align*}
\partial_t U &= \nu\Delta U-\operatorname{div}(U\nabla V),\\
U(0,x) &= \bar m(x)+s n_0(x),\quad x\in\Td.
\end{align*}
Under periodic boundary conditions, this linear parabolic initial-value problem has a unique classical solution $U$.
Thus $\tilde n=U-\bar m$ is uniquely determined.
Finally, given the input $\eta$, solve
\[
\begin{cases}
 \partial_t\tilde q
 =
 \mu\Delta\tilde q
 +2\mu\dfrac{\nabla\bar\theta}{\bar\theta}\cdot\nabla\tilde q
 -\bar\theta\,\tilde q(1+\tilde q)
 -\varepsilon s(1+\tilde q)(\rho*\eta),
\\
 \tilde q(0,x)=s\,q_0(x).
\end{cases}
\]
Let $\tilde q$ be its unique classical solution
\cite{LLMP2005,Friedman1964}.

We next show that a fixed point with $s=1$ yields a solution of (KMR).
Let $(n,q)=\mathcal F((n,q),1)$ be such a fixed point.
By definition, the corresponding $v$ and $(n,q)$ solve the system for the perturbations
with endpoint conditions
$
  v(T)=v_T,\ n(0)=n_0,\ q(0)=q_0
$.
Set
\[
  u=\bar u+\lambda(T-t)+v,\qquad
  m=\bar m+n,\qquad
  \theta=\bar\theta(1+q).
\]
We verify that this change of variables recovers
\eqref{eq:KMR-HJ}--\eqref{eq:KMR-theta} and the prescribed endpoint conditions.

First, consider the equation for $u$.
Using \eqref{eq:sKMR-HJ} and the equation for $v$ at $s=1$, namely
$
-\partial_{t} v-\nu \Delta v-\nabla\bar{u} \cdot \nabla v=\rho \ast (\bar{\theta} q)+\frac{1}{2}|\nabla v|^2
$,
we obtain
\[
\begin{aligned}
 -\partial_t u-\nu\Delta u-\frac12|\nabla u|^2
 &=
 \lambda-\partial_t v-\nu\Delta\bar u-\nu\Delta v
 -\frac12|\nabla\bar u+\nabla v|^2
\\
 &=
 \rho*\bar\theta
 +\left(
 -\partial_t v-\nu\Delta v-\nabla\bar u\cdot\nabla v
 -\frac12|\nabla v|^2
 \right)
\\
 &=
 \rho*\bar\theta+\rho*(\bar\theta q)
 =\rho*\{\bar\theta(1+q)\}
 =\rho*\theta .
\end{aligned}
\]
Thus \eqref{eq:KMR-HJ} holds.
Next, consider the equation for $m$.
The zero-flux identity in Lemma~\ref{lem:fluz} gives
$
  \nu\bar m\nabla{\left(\frac{n}{\bar m}\right)}=\nu\nabla n-n\nabla\bar u
$.
Hence the equation for $n$ becomes
$
  \partial_t n
  =
  \nu\Delta n-\operatorname{div}(n\nabla\bar u)
  -\operatorname{div}(\bar m\nabla v)
  -\operatorname{div}(n\nabla v)
$.
Using \eqref{eq:sKMR-FP}, we obtain
\[
\begin{aligned}
 \partial_t m-\nu\Delta m+\operatorname{div}(m\nabla u)
 &=
 \partial_t n-\nu\Delta n
 +\operatorname{div}(n\nabla\bar u)
 +\operatorname{div}(\bar m\nabla v)
 +\operatorname{div}(n\nabla v)
\\
 &\quad
 -\nu\Delta\bar m+\operatorname{div}(\bar m\nabla\bar u)
 =0 .
\end{aligned}
\]
This proves \eqref{eq:KMR-FP}.
Finally, consider the equation for $\theta$.
Since $\theta=\bar\theta(1+q)$, we have
\[
  \Delta\theta
  =
  \bar\theta\Delta q
  +2\nabla\bar\theta\cdot\nabla q
  +(1+q)\Delta\bar\theta .
\]
Combining this identity with \eqref{eq:sKMR-theta} gives
\[
\begin{aligned}
 \partial_t\theta-\mu\Delta\theta
 &=
 -(1+q)\mu\Delta\bar\theta
 -\bar\theta^2q(1+q)
 {{}-\varepsilon\bar\theta(1+q)(\rho*n)}\\
 &=\bar\theta(1+q)(K-\bar\theta)
 -\bar\theta^2q(1+q)
 -\varepsilon\bar\theta(1+q)\left(\rho*(\bar m+n)\right)
\\
 &=
 \theta(K-\theta)-\varepsilon(\rho*m)\theta .
\end{aligned}
\]
Thus \eqref{eq:KMR-theta} follows.

The definitions $v_T=u_T-\bar u$, $n_0=m_0-\bar m$,
and $q_0=(\theta_0-\bar\theta)/\bar\theta$, together with the change of variables above,
also give $u(T)=u_T$, $m(0)=m_0$, and $\theta(0)=\theta_0$.
Hence $(u,m,\theta)$ constructed from a fixed point with $s=1$
satisfies the original system and the prescribed endpoint conditions.

\begin{proof}[Proof of Proposition~\ref{prop_exsist}]
We first verify that, at $s=0$, for every input
\((\eta,\xi)\), we have
$
  \mathcal{F}((\eta,\xi),0)=(0,0)
$.
Indeed, when \(s=0\), the function \(v\) satisfies
\[
\begin{cases}
-\partial_t v-\nu\Delta v-\nabla\bar u\cdot\nabla v
= \dfrac12|\nabla v|^2,\\
v(T,x)=0.
\end{cases}
\]
The zero function solves this terminal-value problem, so uniqueness gives \(v\equiv0\).
With \(v=0\), the equation for \(\tilde n\) becomes
\[
\begin{cases}
\partial_t\tilde n
= \nu\,\operatorname{div}\left(\bar m\nabla{\left(\frac{\tilde n}{\bar m}\right)}\right),\\
\tilde n(0,x)=0.
\end{cases}
\]
Uniqueness for this linear parabolic initial-value problem gives \(\tilde n\equiv0\).
Finally, the equation for \(\tilde q\) at \(s=0\) is
\[
\begin{cases}
\partial_t\tilde q
= \mu\Delta\tilde q
+2\mu\dfrac{\nabla\bar\theta}{\bar\theta}\cdot\nabla\tilde q
-\bar\theta\,\tilde q(1+\tilde q),\\
\tilde q(0,x)=0.
\end{cases}
\]
The reaction term is locally Lipschitz in \(\tilde q\), so the unique solution with zero initial datum is \(\tilde q\equiv0\).
Thus, for every input \((\eta,\xi)\),
$
  \mathcal{F}((\eta,\xi),0)=(0,0)
$.

{We next exclude fixed points on the boundary of the ball for every $s\in[0,1]$.
Consider a fixed point in the closed ball $\overline{B_\delta^X(0)}$.
Reduce $\bar\delta$ from Section~5.1 if necessary and choose the radius $\delta=\bar\delta$ so that
\[
  2\max\{1,\|\bar m^{-1}\|_{L^\infty}\}\delta
  \le \min\{1/2,1/(2C_R)\}.
\]
Indeed, if $\|(n,q)\|_X\le\delta$, the definition of the norm on $X$ and the bound $\omega_\alpha(t)\le2$ give, for every $t\in[0,T]$,
\[
\begin{aligned}
  \left\|\frac{n(t)}{\bar m}\right\|_{L^\infty}+\|q(t)\|_{L^\infty}
  &\le \max\{1,\|\bar m^{-1}\|_{L^\infty}\}
      \left(\|n(t)\|_{L^\infty}+\|q(t)\|_{L^\infty}\right)\\
  &\le \max\{1,\|\bar m^{-1}\|_{L^\infty}\}\,\omega_\alpha(t)\delta\\
  &\le 2\max\{1,\|\bar m^{-1}\|_{L^\infty}\}\delta
   \le \min\{1/2,1/(2C_R)\}.
\end{aligned}
\]
Thus $r=n/\bar m$ and $q$ satisfy \eqref{hyp:smallness}.}

{Let $(n,q)=\mathcal F((n,q),s)$, and set $\theta=\bar\theta(1+q)$ and $z=\bar\theta q$.
Then $(v,n,\theta)$ satisfies}
  \begin{align}
    -\partial_t v-\nu\Delta v-\nabla\bar u\cdot\nabla v
    &= s\rho*z+\frac12\abs{\nabla v}^2,\\
    \partial_t n
    &= \operatorname{div}\left\{\bar m\left(\nu\nabla{\left(\frac{n}{\bar m}\right)}-\nabla v\right)\right\}
       -\operatorname{div}(n\nabla v),\\
    \partial_t\theta
    &=
    \mu\bar\theta\Delta {q}
    +2\mu\nabla\bar\theta\cdot\nabla {q}
    -\bar\theta^2{q(1+q)}
    -s\eps\bar\theta {(1+q)}(\rho*n)
  \end{align}
{with endpoint data $(v(T),n(0),z(0))=s(v_T,n_0,z_0)$.
This system differs from \eqref{eq:p1}--\eqref{eq:p3} only in that the two nonlocal terms are multiplied by the same factor $s$.
When we differentiate the modified energy, evenness of $\rho$ gives cancellation between
$s\int_{\Td}n(\rho*z)\,dx$ and $-s\int_{\Td}z(\rho*n)\,dx$.
We therefore obtain the same identity as \eqref{eq:energy-identity}.
Since $0\le s\le1$, the norms of the nonlocal terms and endpoint data have the same bounds as at $s=1$.
Thus the improvement argument in Proposition~\ref{prop:n-Linfty-improvement} applies uniformly in $s$.
Choosing the initial and terminal data sufficiently small gives}
\[
\|n(t)\|_{L^\infty}+\|q(t)\|_{L^\infty}
\le \frac{2\delta}{3} (e^{-\alpha t}+e^{-\alpha (T-t)}).
\]
Hence $\|(n,q)\|_X\le 2\delta/3<\delta$, so the fixed point lies in the interior of the ball.

It remains to verify continuity and compactness of the homotopy map.
We first prove continuity of ${\mathcal F}$.
Suppose that, as ${l} \to \infty$,
\[
 \sup_{t\in[0,T]}\|\eta_{l}(t)-\eta(t)\|_{L^\infty}
 +\sup_{t\in[0,T]}\|\xi_{l}(t)-\xi(t)\|_{L^\infty}
 +|s_{l}-s|
 \to 0
\]
and set
$
  (\tilde n_{l},\tilde q_{l})
  :=\mathcal F((\eta_{l},\xi_{l}),s_{l}),
$
$
  (\tilde n,\tilde q)
  :=\mathcal F((\eta,\xi),s)
$.
We must show that, as ${l} \to \infty$,
\[
  \|\tilde n_{l}-\tilde n\|_{L^\infty([0,T]\times\Td)}
  +\|\tilde q_{l}-\tilde q\|_{L^\infty([0,T]\times\Td)}
  \to 0.
\]
We begin with continuous dependence of \(\nabla v_{l}\),
which appears in the equation for \(\tilde n_{l}\).

Let \(v_{l},v\) be the classical solutions of the corresponding Hamilton--Jacobi equations,
\[
\begin{alignedat}{2}
 -\partial_t v_{l}-\nu\Delta v_{l}
 -\nabla\bar u\cdot\nabla v_{l}
 &=s_{l}\rho*(\bar\theta\xi_{l})(t,x)+\frac12|\nabla v_{l}|^2,\qquad
 &v_{l}(T)&=s_{l} v_T,\\
 -\partial_t v-\nu\Delta v-\nabla\bar u\cdot\nabla v
 &=s\rho*(\bar\theta\xi)(t,x)+\frac12|\nabla v|^2,\qquad
 &v(T)&=s v_T .
\end{alignedat}
\]
Reverse time by setting
\[
  \widehat v_{l}(\tau,x):=v_{l}(T-\tau,x),\qquad
  \widehat v(\tau,x):=v(T-\tau,x)
\]
and denote the difference of the forcing terms by
\[
  g_{l}(\tau,x):={ \left[\rho*(\bar\theta(s_l\xi_l-s\xi))\right](T-\tau,x)}.
\]
The difference equation is
\[
\begin{aligned}
 (\partial_\tau-\nu\Delta){(\widehat v_l-\widehat v)}
 -(\nabla\bar u+\nabla\widehat v)\cdot{\nabla(\widehat v_l-\widehat v)}&=g_{l}+\tfrac12{|\nabla(\widehat v_l-\widehat v)|^2},\\
 \qquad {(\widehat v_l-\widehat v)(0)}&=(s_{l}-s)v_T.
\end{aligned}
\]
Differentiating in space and moving the derivative of the quadratic term to the left gives
\[
\begin{aligned}
 &\left(\partial_\tau-\nu\Delta
 -(\nabla\bar u+\nabla\widehat v_{l})\cdot\nabla\right)
 {\nabla(\widehat v_l-\widehat v)}\\
 &\qquad=(D^2\bar u+D^2\widehat v){\nabla(\widehat v_l-\widehat v)}+\nabla g_{l} .
\end{aligned}
\]
Taking the inner product with $\nabla(\widehat v_l-\widehat v)$ and applying Young's inequality, with
$
  C_v:=2\|D^2\bar u+D^2\widehat v\|_{L^\infty([0,T]\times\Td)}+1
$,
we obtain
\[
\begin{aligned}
 &\left(\partial_\tau-\nu\Delta
 -(\nabla\bar u+\nabla\widehat v_{l})\cdot\nabla\right)
 {\left|\nabla(\widehat v_l-\widehat v)\right|^2}
 \leq C_v{\left|\nabla(\widehat v_l-\widehat v)\right|^2}
 +{|\nabla g_l|^2} .
\end{aligned}
\]
The maximum principle gives
\[
 \sup_{\tau\in[0,T]}
 \|\nabla\widehat v_{l}(\tau)-\nabla\widehat v(\tau)\|_{L^\infty}^2
 \leq
 e^{C_vT}
 \left(
 |s_{l}-s|^2\|\nabla v_T\|_{L^\infty}^2
 +T\|\nabla g_{l}\|_{L^\infty([0,T]\times\Td)}^2
 \right).
\]
Moreover,
\[
\begin{aligned}
 {\|\nabla g_l\|_{L^\infty([0,T]\times\Td)}}
 &\leq
 s_{l}\|\nabla\rho\|_{L^1}
 \|\bar\theta\|_{L^\infty}
 \|\xi_{l}-\xi\|_{L^\infty}
 +|s_{l}-s|\|\nabla\rho\|_{L^1}
 \|\bar\theta\xi\|_{L^\infty}
\end{aligned}
\]
and the right-hand side tends to $0$ by uniform convergence of the input sequence and $s_l\to s$.
Time reversal preserves the uniform norm, so the gradient estimate above yields
$
  \|\nabla v_{l}-\nabla v\|_{L^\infty([0,T]\times\Td)}
   \to0
$.

\begingroup
Next, we prove continuous dependence of \(\tilde n_{l}\).
Using the zero-flux identity \(\nu\nabla\bar m=\bar m\nabla\bar u\)
and subtracting the two equations gives
\[
\begin{aligned}
 {(\partial_t-\nu\Delta)(\tilde n_l-\tilde n)}
 &={-\operatorname{div}\left((\tilde n_l-\tilde n)
                    (\nabla\bar u+\nabla v_l)\right)}{-\operatorname{div}\left((\bar m+\tilde n)
                    (\nabla v_l-\nabla v)\right)},\\
 {(\tilde n_l-\tilde n)(0)}&=(s_{l}-s)n_0 .
\end{aligned}
\]
The uniform convergence of \(\nabla v_{l}\) proved above implies that
\(\nabla\bar u+\nabla v_l\) is uniformly bounded in \(l\).
Boundedness of \(\tilde n\) also gives a bound on \(\bar m+\tilde n\).
The periodic heat kernel has \(L^1\) norm \(1\), and its spatial gradient at time \(t>0\)
has \(L^1\) norm at most \(Ct^{-1/2}\).
Let \(e^{\nu t\Delta}\) denote the heat semigroup on the torus.
For any \(0\leq a\leq t\leq T\), Duhamel's formula from time \(a\) gives
\[
\begin{aligned}
 (\tilde n_l-\tilde n)(t)
 &=e^{\nu(t-a)\Delta}(\tilde n_l-\tilde n)(a)\\
 &\quad-\int_a^t \operatorname{div} e^{\nu(t-\tau)\Delta}
   \left[(\tilde n_l-\tilde n)(\nabla\bar u+\nabla v_l)\right](\tau)\,d\tau\\
 &\quad-\int_a^t \operatorname{div} e^{\nu(t-\tau)\Delta}
   \left[(\bar m+\tilde n)(\nabla v_l-\nabla v)\right](\tau)\,d\tau.
\end{aligned}
\]
Applying the $L^\infty$ contraction property of the heat semigroup to the first term
and the heat kernel gradient estimate to the two integrals, we obtain, with a constant \(C>0\) independent of \(l\),
\[
\begin{aligned}
 {\|\tilde n_l(t)-\tilde n(t)\|_{L^\infty}}
 &\leq {\|\tilde n_l(a)-\tilde n(a)\|_{L^\infty}}\\
 &\quad+C\int_a^t(t-\tau)^{-1/2}
 {\|\tilde n_l(\tau)-\tilde n(\tau)\|_{L^\infty}\,d\tau}\\
 &\quad{+C\int_a^t(t-\tau)^{-1/2}
 \|\nabla v_l(\tau)-\nabla v(\tau)\|_{L^\infty}\,d\tau}.
\end{aligned}
\]
Fix \(h\in(0,T]\) such that \(2C\sqrt h\leq 1/2\).
For \(a\leq t\leq\min\{a+h,T\}\), the first integral on the right is bounded by
\[
  2C\sqrt h\sup_{a\leq r\leq\min\{a+h,T\}}\|\tilde n_l(r)-\tilde n(r)\|_{L^\infty}.
\]
Taking the supremum over this interval and absorbing this term into the left-hand side gives
\[
\begin{aligned}
 {\sup_{a\leq t\leq\min\{a+h,T\}}\|\tilde n_l(t)-\tilde n(t)\|_{L^\infty}}
 &\leq {2\|\tilde n_l(a)-\tilde n(a)\|_{L^\infty}}\\
 &\quad+{4C\sqrt h\,\|\nabla v_l-\nabla v\|_{L^\infty([0,T]\times\Td)}} .
\end{aligned}
\]
Iterate this estimate finitely many times from \(a=0\) over intervals of length at most \(h\).
The initial condition and uniform convergence of \(\nabla v_l\) give, for a constant \(C_T>0\) independent of \({l}\),
\[
\begin{aligned}
 \|\tilde n_{l}-\tilde n\|_{L^\infty([0,T]\times\Td)}
 &\leq C_T\left(
 |s_{l}-s|\|n_0\|_{L^\infty}\right.
 \left.+{\|\nabla v_l-\nabla v\|_{L^\infty([0,T]\times\Td)}}
 \right)
\end{aligned}
\]
and the right hand side goes to $0$ as $l \to \infty$.
Here \(T\) is fixed, so dependence of \(C_T\) on \(T\) is allowed.

\par
\endgroup
Finally, we prove continuous dependence of \(\tilde q_{l}\).
The difference equation is
\[
\begin{aligned}
 \partial_t(\tilde q_{l}-\tilde q)
 &=\mu\Delta(\tilde q_{l}-\tilde q)
 +2\mu\frac{\nabla\bar\theta}{\bar\theta}
 \cdot\nabla(\tilde q_{l}-\tilde q)
 \\
 &\quad-\left\{\bar\theta(1+\tilde q_{l}+\tilde q)
 +\varepsilon s_{l}(\rho*\eta_{l})\right\}
 (\tilde q_{l}-\tilde q)
 \\
 &\quad-\varepsilon(1+\tilde q)
 \left\{s_{l}(\rho*\eta_{l})-s(\rho*\eta)\right\},
 \\
 (\tilde q_{l}-\tilde q)(0)&=(s_{l}-s)q_0 .
\end{aligned}
\]
The comparison principle bounds \(\tilde q_{l}\) in \(L^\infty([0,T]\times\Td)\)
uniformly in \({l}\). Since \(\eta_{l}\) is also uniformly bounded, we have
\[
 \sup_{l}
 \|\bar\theta(1+\tilde q_{l}+\tilde q)
 +\varepsilon s_{l}(\rho*\eta_{l})\|_{L^\infty([0,T]\times\Td)}
 <\infty .
\]
Thus the maximum principle gives, for some \(C>0\),
\[
 \|\tilde q_{l}(t)-\tilde q(t)\|_{L^\infty}
 \leq e^{Ct}\left(
 |s_{l}-s|\|q_0\|_{L^\infty}
 +\varepsilon\int_0^t e^{-C\tau}
 \|(1+\tilde q)
 \{s_{l}(\rho*\eta_{l})-s(\rho*\eta)\}(\tau)\|_{L^\infty}
 \,d\tau\right).
\]
Moreover, 
\[
 \|s_{l}(\rho*\eta_{l})-s(\rho*\eta)\|_{L^\infty}
 \leq
 |s_{l}-s|\|\rho\|_{L^1}\|\eta_{l}\|_{L^\infty}
 +s\|\rho\|_{L^1}\|\eta_{l}-\eta\|_{L^\infty}
\]
and the right hand side goes to $0$ as $l \to \infty$.
Thus, letting $l \to \infty$, it follows that
\[
 \|\tilde q_{l}-\tilde q\|_{L^\infty([0,T]\times\Td)}
 {\to0}.
\]
For fixed \(T\), we have \(e^{-\alpha t}+e^{-\alpha(T-t)}\geq2e^{-\alpha T/2}\).
Hence these convergences imply convergence in \(X\), and \(\mathcal F\) is continuous.

Next, we prove compactness of \(\mathcal F\).
Let \(((\eta_{l},\xi_{l}),s_{l})\) be a bounded sequence in
\({\overline{B_\delta^X(0)}}\times[0,1]\), and set
$
 (\tilde n_{l},\tilde q_{l})
 :=\mathcal F((\eta_{l},\xi_{l}),s_{l})
$.
In the Bernstein-type estimate above, take the comparison solution to be \(v=0\) at \(s=0\).
The input sequence is uniformly bounded and \(s_l\in[0,1]\), so the spatial gradients
\(s_l(\nabla\rho)*(\bar\theta\xi_l)\) of the forcing terms are uniformly bounded. Consequently,
$
 \sup_{l}
 \|\nabla v_{l}\|_{L^\infty([0,T]\times\Td)}<\infty
$.

First, the equation for \(\tilde n_{l}\) can be written as
\[
 \partial_t\tilde n_{l}-\nu\Delta\tilde n_{l}
 =-\operatorname{div}\left\{
 \tilde n_{l}(\nu\nabla\log\bar m+\nabla v_{l})
 \right\}
 -\operatorname{div}(\bar m\nabla v_{l}).
\]
Both \(\nu\nabla\log\bar m+\nabla v_l\) and \(\bar m\nabla v_l\) are uniformly bounded.
Using the initial condition \(\tilde n_l(0)=s_l n_0\) and the periodic heat kernel estimates
from the proof of continuous dependence of $\tilde n_l$, we obtain, for a constant \(C_T>0\) independent of \(l\),
\[
 \sup_{l}\|\tilde n_{l}\|_{L^\infty([0,T]\times\Td)}
 \leq {C_T}\left(
 \|n_0\|_{L^\infty}
 +\|\bar m\|_{{L^\infty}}
 \sup_{l}\|\nabla v_{l}\|_{L^\infty([0,T]\times\Td)}
 \right)<\infty .
\]
Thus the vector fields in both divergence terms on the right-hand side of the \(\tilde n_l\) equation
are uniformly bounded in \(l\). The initial data \(s_l n_0\) are also uniformly bounded in \(C^2(\Td)\).
The standard parabolic H\"older estimate \cite[Chapter~III, Theorem~10.1]{LSU1968}
therefore gives some \(\beta\in(0,1)\) such that
$
 \sup_{l}
 \|\tilde n_{l}\|_{C^{\beta/2,\beta}([0,T]\times\Td)}
 {<}\infty
$.

Finally, for \(\tilde q_{l}\), set
$
 \Theta_{l}:=\bar\theta(1+\tilde q_{l})
$.
Then
\[
\begin{cases}
 \partial_t\Theta_{l}-\mu\Delta\Theta_{l}
 =\Theta_{l}\left\{
 K-\varepsilon\rho*(\bar m+s_{l}\eta_{l})-\Theta_{l}
 \right\},
 \\
 \Theta_{l}(0)=(1-s_{l})\bar\theta+s_{l}\theta_0 .
\end{cases}
\]
The comparison principle gives \(\sup_{l}\|\Theta_{l}\|_{L^\infty}<\infty\), and
\[
\begin{aligned}
 &\|K-\varepsilon\rho*(\bar m+s_{l}\eta_{l})-\Theta_{l}\|_{L^\infty}
 \leq
 \|K\|_{L^\infty}
 +\varepsilon\|\rho\|_{L^1}
 \left(\|\bar m\|_{L^\infty}
 +\|\eta_{l}\|_{L^\infty}\right)
 +\|\Theta_{l}\|_{L^\infty}.
\end{aligned}
\]
Each term on the right is uniformly bounded in \(l\).
The same parabolic H\"older estimate therefore gives some $\gamma \in (0,1)$ such that
$
 \sup_{l}
 \|\Theta_{l}\|_{C^{\gamma/2,\gamma}([0,T]\times\Td)}<\infty
$.
Since \(1/\bar\theta\) is smooth, we also obtain
$
 \sup_{l}
 \|\tilde q_{l}\|_{C^{\gamma/2,\gamma}([0,T]\times\Td)}<\infty .
$
The Arzel\`a--Ascoli theorem now implies that
\(\{(\tilde n_{l},\tilde q_{l})\}_{l}\) is relatively compact in \(X\).
Thus \(\mathcal F\) is compact.

{At \(s=0\), the map \(\mathcal F(\cdot,0)\) is zero,
and there are no fixed points on \(\partial B_\delta^X(0)\) for any \(s\in[0,1]\).
The Leray--Schauder fixed point theorem therefore gives a fixed point of \(\mathcal F(\cdot,1)\) in \(B_\delta^X(0)\).
The change of variables above yields a classical solution \((u,m,\theta)\) of (KMR)
with the prescribed endpoint conditions. The choice of the radius ensures \eqref{hyp:smallness},
so the estimates in Proposition~\ref{prop:L2-turnpike_main} also follow.}
\end{proof}
\section*{Acknowledgements}
\noindent
Mamoru Aihara is supported by JST SPRING, Grant Number JPMJSP2119.\\
Shimpei Makida is supported by JSPS KAKENHI Grant Number JP26K24501.

\section*{Use of AI}
The authors used OpenAI Codex and Anthropic Claude Code to translate the manuscript into English and to improve its language and readability.

\bibliographystyle{elsarticle-num}
\bibliography{references}
\end{document}